\documentclass[12pt,twoside,reqno]{amsart}

\usepackage{amsmath,amssymb,amsthm}
\usepackage{fancyhdr}
\usepackage{epsfig,color,colordvi} 

\theoremstyle{plain}
\newtheorem{theorem}{\sc Theorem}[section]
\newtheorem{lemma}{\sc Lemma}[section]
\newtheorem{proposition}[theorem]{\sc Proposition}
\newtheorem{corollary}[theorem]{\sc Corollary}

\theoremstyle{definition}

\newtheorem{example}{\sc Example}[section]

\theoremstyle{remark}
\newtheorem{remark}{\sc Remark}[section]

\newcommand{\be}{\begin{equation}}
\newcommand{\ee}{\end{equation}}
\newcommand{\nn}{\nonumber}

\def\abar{{\bar{a}}}
\def\Gbar{{\bar{G}}}

\def\ombar{{\bar{\om}}}

\def\e{\varepsilon}
\def\ind{\mathbf{1}}  

\DeclareMathOperator{\Leb}{Leb}
\DeclareMathOperator{\sign}{sign}

\def\Fc{\mathcal{F}}  

\def\Hc{\mathcal{H}}
\def\Uc{\mathcal{U}}
\def\Oc{\mathcal{O}}
\def\Pc{\mathcal{P}}
\def\Nc{\mathcal{N}}

\def\atil{{\tilde{a}}}

\def\ytil{{\tilde{y}}}

\def\xtil{{\tilde{x}}}

\def\Xtil{{\widetilde{X}}}

\def\Xbar{{\bar{X}}}
\def\Zbar{{\bar{Z}}}

\def\Sk{{\mathfrak S}}
\def\Hvec{\mathbf{H}}
\def\Yvec{\mathbf{Y}}

\newcommand{\Ev}{{\bf E}}            

\newcommand{\Pv}{{\bf P}}

\newcommand{\Pf}{{\mathbb P}}
\newcommand{\Ef}{{\mathbb E}}
\newcommand{\Vv}{{\text{\bf Var}}}

\newcommand{\Zb}{\mathbb Z}
\newcommand{\Rb}{\mathbb R}
\newcommand{\Nb}{\mathbb N}

\newcommand{\si}{\sigma}
\newcommand{\om}{\omega}
\newcommand{\Om}{\Omega}

\newcommand{\te}{\theta}
\newcommand{\de}{\delta}

\newcommand{\al}{\alpha}
\newcommand{\ga}{\gamma}
\newcommand{\vr}{\varrho}
\newcommand{\ve}{\varepsilon}

\newcommand{\tfl}[1]{\lfloor{#1}\rfloor}   

\newcommand{\wt}{\widetilde}

\def\bbR{\mathbb R}

\def\bbQ{\mathbb Q}

\def\fQ{{{\rm Q}\kern-.65em {}^{{}_/ }\,}}
\def\fQQ{ {{\rm Q}\kern-.57em \scriptscriptstyle{}^{]\kern.055em[}\,}}

\def\w{\omega}

\def\ord{\kern0.1em o\kern-0.02em{}_{\ds\breve{}}\kern0.1em}
\def\Ord{\kern0.1em O\kern-0.02em{\ds\breve{}}\kern0.1em}
\def\ds{\displaystyle}

\def\fmonth{\ifcase\month\or Jan\or Feb\or Mar\or Apr
\or May\or Jun\or Jul\or Aug\or Sep
\or Oct\or Nov\or Dec\fi\ }
\def\mmddyyyy{\the\month.\the\day.\the\year}
\def\ddmmyyyy{\the\day.\the\month.\the\year}
\def\Mddyyyy{\fmonth~\the\day,~\the\year}

\providecommand{\abs}[1]{\left\vert#1\right\vert}
\providecommand{\norm}[1]{\left\Vert#1\right\Vert}

\numberwithin{equation}{section}

\allowdisplaybreaks[1]

\begin{document}

\author{M\'arton Bal\'azs} 
\address{M.~Bal\'azs, Mathematics Department, University 
of Wisconsin-Madison,
Van Vleck Hall, Madison, WI 53706, USA.}
\thanks{M.~Bal\'azs was partially supported by 
Hungarian Scientific Research Fund (OTKA) grant
T037685.}
\email{balazs@math.wisc.edu}
\author{Firas Rassoul-Agha}
\address{F. Rassoul-Agha, Mathematical Biosciences Institute, 
Ohio State University,
231 West 18$^{\hbox{th}}$ Avenue, Columbus, OH 43210, USA.}
\email{firas@math.ohio-state.edu}
\author{Timo Sepp\"al\"ainen}
\address{T. Sepp\"al\"ainen, Mathematics Department, University 
of Wisconsin-Madison,
Van Vleck Hall, Madison, WI 53706, USA.}
\email{seppalai@math.wisc.edu}
\thanks{T.~Sepp\"al\"ainen was partially supported by
National Science Foundation grant DMS-0402231.}

\date{Submitted September 7, 2005. Revised January 11, 2006.}


\begin{abstract}
We study space-time fluctuations around a characteristic line for a
one-dimensional interacting system  known as the random
average process. The state of this system is a real-valued 
function on the integers.  New values of the function are
created by averaging previous values with random weights. 
  The fluctuations analyzed occur on the scale
$n^{1/4}$ where $n$ is the ratio of macroscopic and microscopic
scales in the system. The  limits of the fluctuations
 are described 
by a family of 
Gaussian processes. In  cases of known product-form 
invariant distributions, this limit is a two-parameter process
whose time marginals are fractional Brownian
motions with Hurst parameter $\tfrac14$. 
Along the way we study the limits of quenched mean processes
for a random walk in a space-time random environment. 
These limits also happen at scale $n^{1/4}$ and are described
by certain Gaussian processes that we identify. In particular, 
when we look at a backward quenched mean process, the limit 
process is the solution  of a stochastic heat equation. 
\end{abstract}




\title[Random average process]{The  random average
process and  random walk 
in a space-time random environment in one dimension }

\maketitle



\section{Introduction}

\subsection*{Fluctuations for asymmetric interacting systems}
An asymmetric interacting system is a random process  
$\sigma_\tau=\{\sigma_\tau(k): k\in\mathcal{K}\}$
of many components $\sigma_\tau(k)$ 
 that influence each others' evolution.
 Asymmetry means here that the components have 
 an average drift in some spatial direction.  Such processes
are called interacting particle systems because often these
 components can be thought of  as particles. 
   
To orient the reader, let us first think of 
  a single  random walk 
$\{X_\tau: \tau=0,1,2,\dotsc\}$
that evolves by itself.  
  For random walk we scale both space and 
time by $n$ because on this scale we see the long-term velocity:
$n^{-1}X_{\tfl{nt}}\to tv$ as $n\to\infty$ where $v=EX_1$.  
 The random walk is {\sl diffusive}   which  means that its fluctuations
occur on the scale $n^{1/2}$, as revealed by the 
classical central limit theorem: 
$n^{-1/2}(X_{\tfl{nt}}-ntv)$ converges weakly to a Gaussian distribution. 
The Gaussian limit is {\sl universal} here because it arises 
regardless of the choice of step distribution for the random
walk, as long as a  square-integrability  hypothesis is satisfied. 

For asymmetric interacting systems we typically also scale time and 
space by the same factor $n$, and this is known as 
{\sl   Euler scaling.}     However, in certain classes of 
 one-dimensional asymmetric interacting systems the random
evolution produces 
fluctuations   of smaller
order than the natural diffusive scale. 
Two types of such phenomena have been discovered. 

\medskip

(i) 
In Hammersley's process, in asymmetric exclusion, and 
in some other closely related systems,
dynamical fluctuations occur on the scale $n^{1/3}$. 
 Currently known rigorous 
results
 suggest  that the Tracy-Widom 
distributions from random matrix theory are the universal
limits of these $n^{1/3}$ fluctuations.  

The seminal 
works in this context are by Baik, Deift and Johansson
\cite{baik-deif-joha} on Hammersley's
process and by Johansson \cite{joha} on the  exclusion process.  
We should point out though that \cite{baik-deif-joha} does
not explicitly discuss Hammersley's process, but instead
the maximal number of planar Poisson points on an increasing
path in a rectangle.  One
can intrepret the results in \cite{baik-deif-joha} as
fluctuation results for Hammersley's process with a 
special initial configuration. The connection between 
the increasing path model and Hammersley's process goes
back to Hammersley's paper \cite{hamm}. It was first utilized
 by Aldous and Diaconis \cite{aldo-diac95} (who also
named the process), and then further in the papers
\cite{sepp96, sepp98ldp}. 

\medskip

(ii) The second type has fluctuations of the 
order $n^{1/4}$ and limits described by a family of
self-similar  Gaussian processes
that includes fractional Brownian motion with Hurst parameter
$\tfrac14$.   This result was first proved for a system of 
independent random walks \cite{sepp-rw}. One of the main results
of the current paper shows that the $n^{1/4}$ fluctuations
also appear in a family of interacting systems called 
{\sl random average processes} in one dimension. 
The same family of  limiting Gaussian processes appears
here too, suggesting that these
limits are universal for some class of interacting systems. 

\medskip

The random average processes (RAP) studied in the present paper
describe a random   real-valued 
 function on the integers whose values evolve by 
jumping to random convex combinations of  values in a finite 
neighborhood. It could be thought of as a caricature model 
for an interface between two phases on the plane, hence
we call the state a {\sl height function}. 
RAP is  
  related to the so-called linear systems 
discussed in Chapter IX of Liggett's monograph \cite{ligg85}. 
RAP was introduced by Ferrari and Fontes \cite{ferr-font-rap}
who studied  the fluctuations from initial linear slopes. 
In particular, they discovered that the height over the origin
satisfies a  central limit theorem in the time scale $t^{1/4}$.
The Ferrari-Fontes results suggested 
RAP to us as a fruitful place to investigate whether 
the $n^{1/4}$ fluctuation picture discovered in \cite{sepp-rw}
for independent
walks had any claim to universality.

There are two ways to see the lower order dynamical fluctuations.
\begin{itemize}
\item[{(1)}]  One can take deterministic 
initial conditions so that only dynamical randomness 
is present.  \item[{(2)}] Even if the initial state is random with central
limit scale fluctuations,  one can find the lower order fluctuations
 by looking  at the evolution of the process 
 along a  characteristic curve. 
\end{itemize} 
Articles \cite{baik-deif-joha} and   
  \cite{joha} studied the evolutions of special deterministic
initial states of 
Hammersley's process  and  the exclusion process. 
Recently Ferrari and Spohn \cite{ferr-spoh} have extended this
analysis to the fluctuations across a characteristic
in a stationary exclusion process.  
 The 
general nonequilibrium hydrodynamic limit situation is
still out of reach for these models.  \cite{sepp-rw} contains
a tail bound for Hammersley's process
 that suggests $n^{1/3}$ scaling also in the nonequilibrium
situation, including along a  shock which can be regarded as  a
``generalized'' characteristic.  
 
Our results for the random average process are for 
the general hydrodynamic limit setting. 
The initial increments
of the random height function are assumed
independent 
and subject to some moment bounds. Their means and variances must vary
 sufficiently regularly   
to satisfy a H\"older condition.    Deterministic initial 
increments qualify here as a special case of independent.

The classification of the systems mentioned above
(Hammersley, exclusion, independent walks, RAP)
 into $n^{1/3}$ and $n^{1/4}$ fluctuations
coincides with their classification according to 
  type of macroscopic equation. 
Independent particles and RAP 
are macroscopically governed by linear first-order partial differential
equations  $u_t+bu_x=0$.  In contrast, macroscopic evolutions
of Hammersley's process and the 
exclusion process obey genuinely 
nonlinear Hamilton-Jacobi equations $u_t+f(u_x)=0$
that create shocks. 
 
Suppose we start off one of these systems so that the initial
state fluctuates
on the $n^{1/2}$ spatial scale, for example in an stationary distribution. 
Then   the fluctuations
of the entire system    on the    $n^{1/2}$  scale  
simply consist of initial fluctuations transported
along the deterministic characteristics of the macroscopic equation. 
This is a consequence of the lower order  of dynamical fluctuations. 
When the macroscopic equation is linear this is the whole picture of
diffusive fluctuations. 
In the nonlinear case the behavior at the shocks (where 
characteristics merge) also needs to be resolved. This has been
done for  the exclusion
process   \cite{reza02tasep} and for Hammersley's process
  \cite{sepp02diff}.

\subsection*{Random walk in a space-time random
environment} 
Analysis of the random average 
process  utilizes a dual description in terms of 
backward random walks in a space-time random
environment.  
Investigation of the fluctuations of RAP  leads to a study of fluctuations
of these random walks,  both quenched invariance principles for 
the walk itself and limits for the quenched mean process. 
The quenched  invariance principles have been reported
elsewhere \cite{qclt-spacetime}.  The results for the 
quenched mean process are included in the present paper because
they are intimately connected to the random average process
results. 

We look at two types of processes of quenched means. We call them
forward and backward.  In the forward case the initial point of
the walk is fixed, and the walk runs for a specified amount of time
 on the space-time lattice.  In the backward case the initial
point moves along a characteristic, and the walk runs until it
reaches the horizontal axis. Furthermore, in both cases we let
the starting point vary horizontally (spatially), and so we 
have  a space-time process.  In both cases we describe a limiting
Gaussian process, when space is scaled by $n^{1/2}$, time
by $n$, and the magnitude of the fluctuations by $n^{1/4}$. 
In particular, in the backward case we find a limit process 
that solves the stochastic heat equation. 

There are two earlier papers on the 
quenched mean of this
random walk in a space-time random environment. These previous
results  were proved under 
assumptions of small enough noise and finitely many possible
values for the random probabilities. 
Bernabei \cite{bern} showed that the centered quenched mean,
normalized by its own standard deviation, converges
to a normal variable. Then  separately he showed that 
this standard deviation is bounded above and below 
on the order $n^{1/4}$.  Bernabei has results also in 
dimension 2, and also for the quenched covariance of the walk.
Boldrighini and Pellegrinotti \cite{bold-pell-01}
also proved a normal limit in the scale $n^{1/4}$ 
for what they term the ``correction'' caused by the 
random  environment on the mean of 
a test function.

\subsection*{Finite-dimensional versus process-level convergence}
Our main results all state that the finite-dimensional 
distributions of a process of interest converge to the 
finite-dimensional distributions of a certain Gaussian 
process specified by its covariance function. 
We have not proved process-level tightness, except in the case
of forward quenched means for the random walks 
where we compute a bound on the sixth moment of the 
process increment.

\subsection*{Further relevant literature}
It is not clear what exactly are the systems ``closely related''
to Hammersley's process or exclusion process, alluded to
in the beginning of the Introduction,  that share the
$n^{1/3}$ fluctuations and Tracy-Widom limits.  The processes
for which rigorous proofs exist all have an underlying
representation in terms of a last-passage percolation model.
Another such example is ``oriented digital boiling''  
studied by Gravner, Tracy and Widom \cite{grav-trac-wido}.
(This model was studied earlier in \cite{sepp98perc} 
and \cite{joha01}  under different names.) 

Fluctuations of the current were initially studied from 
the perspective of a moving observer traveling with a general
speed. The fluctuations are diffusive, and the limiting
variance is a function of the speed of the observer.
 The special nature of the characteristic speed 
manifests itself in the vanishing of the limiting variance
on this diffusive scale.
The early paper of 
Ferrari and Fontes \cite{ferr-font-94} 
 treated the asymmetric exclusion process.
Their work was extended  by Bal\'azs \cite{bala-03} 
 to a class of deposition models  that includes the 
much-studied zero range process and a generalization  
called the  bricklayers' process. 

Work on the fluctuations of  Hammersley's process 
and the exclusion process has connections to several 
parts of mathematics. Overviews of some of these 
links appear in  papers \cite{aldo-diac99}, \cite{deif},
and  \cite{groe}.  
General treatments of large scale behavior of interacting
random systems can be found in \cite{dema-pres}, \cite{kipn-land}, 
\cite{ligg85}, \cite{ligg99},
 \cite{spoh}, and \cite{vara}.

\subsection*{Organization of the paper} 
We begin with the  description of the random average process
and the limit theorem for it in Section \ref{rap-section}.
Section \ref{rwre-section} describes the 
  random walk in a space-time random environment and 
the limit theorems for quenched mean processes.  
The proofs begin with Section \ref{prelim-sec} that lays out
some preliminary facts on random walks.  
Sections \ref{bw-walk-pf-sec} and \ref{fw-walk-pf-sec}
prove the fluctuation results for random walk, and the final Section
\ref{pf-rap-sec}  proves the limit theorem for RAP. 

The reader only interested in the 
 random walk can read Section \ref{rwre-section} and the 
proofs for the random walk limits   independently of the rest 
of the paper, except for
certain definitions and a hypothesis  which have been labeled. 
The RAP results can be read independently of the 
random walk, but their proofs  depend on the random walk results. 

\subsection*{Notation} We summarize here some notation and 
conventions  for 
quick reference. The set of natural numbers is 
 $\Nb=\{1,\,2,\,3,\,\dotsc\}$, while $\Zb_+=\{0,\,1,\,2,\,3,\,\dotsc\}$
and $\Rb_+=[0,\infty)$.
On the two dimensional integer lattice $\Zb^2$ standard  basis 
vectors are $e_1=(1,0)$ and $e_2=(0,1)$.  The $e_2$-direction
 represents time. 

 We need several different probability 
measures and corresponding expectation operators. 
$\Pf$ (with expectation $\Ef$)
 is the probability measure on the space $\Om$ of
environments $\om$. $\Pf$ is an i.i.d.\ product measure
across the coordinates indexed by the  space-time lattice $\Zb^2$.
 $\Pv$  (with expectation $\Ev$) is the 
probability measure of the initial state of the random
average process.  $\Ev^\om$ is used to emphasize that 
an expectation over initial states is taken with a 
fixed environment $\om$.
Jointly the environment and initial state are independent,
so the joint measure is the product $\Pf\otimes\Pv$. 
$P^\om$ (with expectation $E^\om$) is the quenched path
measure of the random walks in environment $\om$. The
annealed measure for the walks is $P=\int P^\om\,\Pf(d\om)$. 
Additionally, we use $P$ and $E$ for generic probability
measures and expectations for processes that are not 
part of this specific set-up, such as Brownian motions 
and limiting Gaussian processes. 
  
The environments $\om\in\Om$ are configurations
$\om=(\om_{x,\tau}: (x,\tau)\in\Zb^2)$  of vectors indexed by the
space-time lattice $\Zb^2$.   Each element $\om_{x,\tau}$ is a 
  probability vector of length 
 $2M+1$, denoted also by $u_\tau(x)=\om_{x,\tau}$, 
and in terms of coordinates $u_\tau(x)=(u_\tau(x,y): -M\leq y\leq M)$.
The environment at a fixed time value $\tau$ is 
$\ombar_\tau=(\om_{x,\tau}: x\in\Zb)$.  Translations on $\Omega$
are defined by $(T_{x,\tau}\om)_{y,s}=\om_{x+y,\tau+s}$. 

 $\tfl{x}=\max\{n\in\Zb: n\leq x\}$ is the lower integer part
of a real $x$.
Throughout, $C$ denotes a constant whose exact  value
is immaterial and  can change from line
to line. 
The density and cumulative distribution function of 
the  centered Gaussian distribution with variance $\sigma^2$
are denoted by $\varphi_{\sigma^2}(x)$ and $\Phi_{\sigma^2}(x)$. 
 $\{B(t):t\geq 0\}$ is one-dimensional
standard Brownian motion, in other words the Gaussian process
with covariance $EB(s)B(t)=s\land t$.  

\subsection*{Acknowledgements} The authors thank P.~Ferrari and 
L.~R.~Fontes for comments on article \cite{ferr-font-rap}
and J.~Swanson for helpful discussions.

\def\Ybar{\bar Y}
\def\ybar{\bar y}

\def\qbar{{\bar{q}}}
\def\ombar{\bar{\omega}}
\def\rhobar{\bar{\rho}}
\def\rg{M}  

\section{The random average process} \label{rap-section}

The state of the  random average process (RAP) 
is a height function  $\si: \Zb\to\Rb$. It can also
be thought of as a sequence  
$\si=(\si(i)\,:\,i\in\Zb)\in\Rb^\Zb$ where $\si(i)$ is
the height of an interface above site $i$. 
The state evolves  in discrete time  according to the following rule. 
At each time
point $\tau=1,2,3,\dotsc$ and 
at each site $k\in\Zb$, a  random  probability vector 
$u_\tau(k)=(u_\tau(k,j): -\rg\leq j\leq \rg )$  of length $2\rg+1$ is drawn. 
Given the state $\si_{\tau-1}=(\si_{\tau-1}(i)\,:\,i\in\Zb)$ 
at time $\tau-1$,  
the height value at site $k$ is then updated to 
 \be
\si_\tau(k)=\sum_{j: \abs{j}\leq \rg}u_\tau(k,\,j) \si_{\tau-1}(k+j).
\label{si-dyn-1}
\ee
This update is performed independently at each site $k$ to form
the state $\si_\tau=(\si_\tau(k)\,:\,k\in\Zb)$ at time $\tau$. 
The same step is repeated at the next time $\tau+1$ with 
new independent draws of the probability vectors. 

So, given an initial state $\si_0$,
the process $\sigma_\tau$  is constructed with a collection 
$\{ u_\tau(k): \tau\in\Nb,\, k\in\Zb\}$  of independent and identically 
distributed  random vectors.  
   These random vectors are defined on  a probability space
$(\Om, \Sk, \Pf)$.  If $\sigma_0$ is also random with 
distribution $\Pv$, then  $\sigma_0$ and the vectors 
$\{u_\tau(k)\}$ are independent, in other words the joint
distribution is $\Pf\otimes\Pv$. 
 We write $u^\om_\tau(k)$ to make 
  explicit  the dependence on $\om\in\Om$.  
$\Ef$ will denote expectation under the measure $\Pf$. 
  $\rg$ is the {\sl range} and is a fixed finite parameter of
the model.  $\Pf$-almost surely each random vector
$u_\tau(k)$ satisfies 
 \[
\text{$0\leq u_\tau(k,j)\leq 1 $ for all $-\rg\leq j\leq \rg$, and }\; 
\sum_{j=-\rg}^{\rg}u_\tau(k,j)=1.  
\]
It is often convenient to allow values $u_\tau(k,j)$ 
for all $j$. Then automatically $u_\tau(k,j)=0$ 
for $\abs{j}>\rg$. 

Let \[p(0,j)=\Ef u_0(0,j)\] denote the averaged probabilities. 
Throughout the paper we make two fundamental assumptions.  

(i) First,    there is no 
integer $h>1$ such that, for some $x\in\Zb$, 
\[
\sum_{k\in\Zb} p(0, x+kh)=1.
\]
This is also expressed by saying that the {\sl span}
 of the random walk with  
jump probabilities $p(0,j)$  is 1 \cite[page 129]{durr}.
It follows that the group generated by $\{x\in\Zb: p(0,x)>0\}$ 
is all of $\Zb$, in other words this walk is {\sl aperiodic} in 
Spitzer's terminology \cite{spit}. 

(ii) Second, we assume that
\be
\Pf\{ \max_j u_0(0,j) < 1\}>0. 
\label{ellipt}
\ee
If this  assumption fails, then $\Pf$-almost surely 
for each $(k,\tau)$ 
there exists $j=j(k,\tau)$ such that 
$u_\tau(k,j)=1$. No averaging happens, but instead 
$\si_\tau(k)$ adopts the value $\si_{\tau-1}(k+j)$. 
The behavior is then different from that described by
our results.

No further hypotheses are required of the distribution $\Pf$
on the probability vectors.  Deterministic weights
$u^\om_\tau(k,j)\equiv p(0,j)$  are also admissible, in which
case \eqref{ellipt} requires $\max_j p(0,j)<1$.

In addition to the height process $\si_\tau$
 we also consider  the increment process
$\eta_\tau=(\eta_\tau(i):i\in\Zb)$ defined by 
\[
\eta_\tau(i)=\si_\tau(i)-\si_\tau(i-1).
\]
From \eqref{si-dyn-1} one can deduce a similar
linear equation for 
the evolution of the increment process.  
However, the weights are not necessarily nonnegative,
and even if they are, they do not necessarily sum to
one.

Next we define several constants that appear in the results.  
\be
D(\om)=\sum_{x\in\Zb} x \; u^\om_0(0,x)
\label{def-D}
\ee
is the drift  at the origin. Its 
mean is $V=\Ef(D)$ and  variance
\be
\sigma_D^2=\Ef[(D-V)^2].
\label{def-si_D}
\ee
A variance under    averaged probabilities    is computed by 
\be
\sigma_a^2=\sum_{x\in\Zb} (x-V)^2\; p(0,x).
\label{def-si_a}
\ee
Define random and averaged characteristic functions by
\be
\phi^\om(t)=\sum_{x\in\Zb} u_0^\om(0,x)e^{itx}
\quad\text{and}\quad 
\phi_a(t)= \Ef \phi^\om(t)= \sum_{x\in\Zb} p(0,x)e^{itx},
\label{def-chf}
\ee
and then further 
\be
\lambda(t)=\Ef[\,\lvert\phi^\om(t)\rvert^2\,]
\quad\text{and}\quad 
\bar{\lambda}(t)=\lvert\phi_a(t)\rvert^2. 
\label{def-chf-2}
\ee
Finally, define a positive constant $\beta$  by
\be
\beta= \frac1{2\pi} \int_{-\pi}^\pi
 \frac{1- \lambda(t)} {1- \bar{\lambda}(t)} \,dt.
\label{def-beta-1}
\ee 
The assumption of span 1 implies that $\lvert\phi_a(t)\rvert=1$
only at multiples of $2\pi$. Hence the integrand above
is positive at $t\neq 0$. Separately one can check that the 
integrand has  a finite limit as $t\to 0$. Thus $\beta$ is
well-defined and finite. 

In Section \ref{prelim-sec}   we can give these  constants,
especially $\beta$,    more probabilistic meaning 
from the perspective of  the underlying
 random walk in random environment. 
 
For the limit theorems we consider 
  a sequence $\si^n_\tau$ of the random average processes,
indexed by $n\in\Nb=\{1,\,2,\,3,\,\dotsc\}$.  
  Initially we set $\si^n_0(0)=0$. 
For each $n$ we assume that the initial increments
$
\left\{\eta_0^n(i): i\in\Zb\right\}
$
are independent random variables, with
\be
\Ev [\eta_0^n(i)]=\vr(i/n)\quad \text{and}
\quad 
\Vv[\eta_0^n(i)]=v(i/n). 
\label{eq:meta}
\ee
The functions  $\vr$ and $v$ that appear above
 are assumed to be   uniformly bounded functions on $\Rb$ and
to satisfy this  
 local H\"older continuity: 
\be
\begin{split} 
&\text{For each 
compact interval $[a,b]\subseteq\Rb$ there exist  } \\
&\text{$C=C(a,b)<\infty$ and  $\gamma=\gamma(a,b)>1/2$  such that }\\
 & \abs{\vr(x)-\vr(y)}+\abs{v(x)-v(y)}\leq C\abs{x-y}^\gamma 
\quad\text{for $x,y\in[a,b]$.}
\end{split}
\label{holder}
\ee
The function $v$ must be nonnegative, but the sign of $\vr$ is not
restricted.  Both functions are allowed to vanish. 
In particular, our hypotheses permit   deterministic initial heights
which implies that $v$ vanishes identically. 

 The distribution on initial heights and increments described above
 is denoted by $\Pv$. 
We make this uniform  moment hypothesis on the increments: 
\be
\text{there exists $\alpha>0$ 
such that } \quad
 \sup_{n\in\Nb,\,i\in\Zb} \Ev[\,\lvert\eta_0^n(i)\rvert^{2+\alpha} \,]
<\infty. 
\label{eta-ass-p}
\ee
We assume that the processes $\si^n_\tau$ are all defined on 
the same probability space. 
 The environments
$\om$ that drive the dynamics are independent of the initial
states $\{\si^n_0\}$, so the joint distribution of $(\om, \{\si^n_0\})$ is 
$\Pf\otimes\Pv$.  When computing an expectation under
a fixed $\om$ we write  
 $\Ev^\om$. 

On the larger space and time scale the height function 
is simply rigidly translated at speed $b=-V$, and the same is 
also true of the central limit fluctuations of the 
initial height function. Precisely speaking,
 define a function $U$ on $\Rb$ by $U(0)=0$ and
$U'(x)=\vr(x)$. Let $(x,t)\in\Rb\times\Rb_+$. 
The assumptions made thus far imply that both 
\be
n^{-1}\si^n_{\tfl{nt}}(\tfl{nx})\longrightarrow U(x-bt)
\label{si-hydro}
\ee
and 
\be
\frac{\si^n_{\tfl{nt}}(\tfl{nx})- nU(x-bt)}{\sqrt{n}}
\,-\, \frac{\si^n_{0}(\tfl{nx}-\tfl{nbt})- nU(x-bt)}{\sqrt{n}}
\label{si-hyd-fl} \longrightarrow 0 
\ee
in probability, as $n\to\infty$. (We will not give a proof.
This follows from easier versions of the estimates in the paper.)  
Limit \eqref{si-hydro}
 is the ``hydrodynamic limit'' of the process. The large scale
evolution of the height process
 is thus governed by the linear transport equation 
\[
w_t+bw_x=0. 
\]
This equation  is uniquely
 solved by $w(x,t)=U(x-bt)$ given the initial function
$w(x,0)=U(x)$. 
The lines $x(t)=x+bt$ are the characteristics 
of this equation, the curves along which
the equation carries information. 
Limit \eqref{si-hyd-fl} says that fluctuations on the diffusive scale 
do not include any randomness from the evolution, only a 
translation of initial  fluctuations along characteristics. 

We find interesting height fluctuations along a
macroscopic characteristic 
line $x(t)=\ybar+bt$, and around such a line 
 on the microscopic spatial scale $\sqrt{n}$. 
The magnitude of these fluctuations is of the order 
$n^{1/4}$,  so we study the process 
\[
z_n(t,r)=n^{-1/4}
\bigl\{\si_{\tfl{nt}}^n(\tfl{n\ybar}+\tfl{r\sqrt{n}\,}+\tfl{ntb})
-\si_0^n(\tfl{n\bar y}+\tfl{r\sqrt{n}\,})\bigr\},
\]
indexed by $(t,r)\in\Rb_+\times\Rb$,  
for a fixed $\bar y\in\Rb$. In terms of the increment 
process $\eta^n_\tau$, 
$z_n(t,0)$ 
is the net flow  from right to left across the discrete 
characteristic $\tfl{n\ybar}+\tfl{nsb}$,
during the time interval $0\leq s\leq t$. 

Next we describe the limit of $z_n$. 
Recall the constants defined in \eqref{def-si_D}, \eqref{def-si_a}, 
and \eqref{def-beta-1}. Combine them into a new  constant
\be
\kappa=\frac{\sigma_D^2}{\beta\sigma_a^2  }.
\label{def-kappa}
\ee
Let  $\{B(t):t\geq 0\}$ be one-dimensional  standard  Brownian motion.
Define two functions $\Gamma_q$ and $\Gamma_0$ on
$(\Rb_+\times\Rb)\times(\Rb_+\times\Rb)$:
\be
\Gamma_q((s,q),(t,r)) = 
\frac\kappa{2}
\int_{ \sigma_a^2\lvert t-s\rvert}
^{\sigma_a^2(t+s)} \frac1{\sqrt{2\pi v}}
\exp\Bigl\{-\frac1{2v}(q-r)^2\Bigr\} \,d v
\label{y-cov}
\ee
and 
\be
\begin{split}
&\Gamma_0((s,q),(t,r))
=\int_{q\lor r}^\infty 
P[ \sigma_a  B(s)> x-q]P[  \sigma_a B(t)> x-r]\,dx \\
&\qquad\qquad 
-  \Bigl\{ \ind_{\{r>q\}} \int_q^r P[   \sigma_a B(s)> x-q]
P[  \sigma_a B(t)\leq x-r]\,dx \\
&\qquad\qquad\qquad\qquad + 
\ind_{\{q>r\}} \int_r^q P[   \sigma_a B(s)\leq x-q]
P[  \sigma_a B(t)> x-r]\,dx \Bigr\} \\
&\qquad\qquad + 
 \int_{-\infty}^{q\land r} P[   \sigma_a B(s)\leq x-q]
P[  \sigma_a B(t)\leq x-r]\,dx. 
\end{split}
\label{def-Ga_0}
\ee
The boundary values are such that 
$\Gamma_q((s,q),(t,r))=\Gamma_0((s,q),(t,r))=0$ if either $s=0$
or $t=0$. 
We will see later that 
$\Gamma_q$ is the limiting covariance of the backward quenched mean process
of a related random walk in random environment. 
$\Gamma_0$ is the covariance for fluctuations contributed by the 
initial increments of the random average process. (Hence the
subscripts $q$ for quenched and $0$ for initial time. The 
subscript on $\Gamma_q$ has nothing to do with the argument $(s,q)$.) 

The integral expressions  above are the form in which 
$\Gamma_q$ and $\Gamma_0$ appear in the 
proofs. For $\Gamma_q$ the key point is the limit 
\eqref{G-q-appear} which is evaluated earlier in  \eqref{G-lim1}.  
$\Gamma_0$ arises in Proposition \ref{Y-findim-prop}. 

Here are alternative succinct representations for 
$\Gamma_q$ and $\Gamma_0$. 
Denote   the centered Gaussian density with variance $\sigma^2$
by 
\be
\varphi_{\sigma^2}(x)=\frac1{\sqrt{2\pi \sigma^2}}
\exp\Bigl\{-\frac1{2\sigma^2}x^2\Bigr\}  
\label{gauss-fi}
\ee 
and its distribution function by
$
\Phi_{\sigma^2}(x)=\int_{-\infty}^x 
\varphi_{\sigma^2}(y) \,dy.
$
Then define
\[
\Psi_{\sigma^2}(x)= \sigma^2 \varphi_{\sigma^2}(x) -x(1-\Phi_{\sigma^2}(x)),
\]
which is an antiderivative of $\Phi_{\sigma^2}(x)-1$.  In these terms, 
\[
\Gamma_q((s,q),(t,r)) 
= \kappa \Psi_{\sigma_a^2(t+s)} \bigl(\lvert q-r\rvert\bigr)
-
\kappa \Psi_{\sigma_a^2\lvert t-s\rvert} \bigl(\lvert q-r\rvert\bigr).
\]
and
\[
\Gamma_0((s,q),(t,r)) 
=  \Psi_{\sigma_a^2s} \bigl(\lvert q-r\rvert\bigr)+
\Psi_{\sigma_a^2t} \bigl(\lvert q-r\rvert\bigr)
-  \Psi_{\sigma_a^2(t+s)} \bigl(\lvert q-r\rvert\bigr).
\]

\begin{theorem} Assume {\rm \eqref{ellipt}}
and that  the averaged probabilities
$p(0,j)=\Ef u^\om_0(0,j)$ have lattice span 1. 
Let $\vr$ and $v$ be two uniformly bounded functions 
on $\Rb$ that satisfy the local H\"older condition 
{\rm \eqref{holder}}.  For each $n$, let $\si^n_\tau$ 
be a random average process normalized by $\si^n_0(0)=0$
and whose initial increments $\{\eta^n_0(i):i\in\Zb\}$ 
are independent and satisfy {\rm\eqref{eq:meta}} and 
{\rm\eqref{eta-ass-p}}. Assume the environments $\om$
independent of the initial heights $\{\si^n_0: n\in\Nb\}$.

Fix $\ybar\in\Rb$. 
Under the above assumptions  
the finite-dimensional distributions 
of the  process $\{z_n(t,r): (t,r)\in\Rb_+\times\Rb \}$ converge weakly 
as $n\to\infty$ 
to the finite-dimensional distributions 
of the   mean zero Gaussian
process $\{z(t,r): (t,r)\in\Rb_+\times\Rb \}$ specified by the covariance 
\be
\begin{split}
Ez(s,q)z(t,r)&= 
 \vr(\ybar)^2 \Gamma_q((s,q),(t,r))
+ v(\ybar) \Gamma_0((s,q),(t,r)).
\end{split}
\label{rap-cov}
\ee
\label{rap-thm-1}
\end{theorem}

The statement means that, given  space-time points 
$(t_1,r_1),\dotsc,(t_k,r_k)$, the $\Rb^k$-valued random vector 
$(z_n(t_1,r_1),\dotsc, z_n(t_k,r_k))$ converges in distribution to 
the random vector $(z(t_1,r_1),\dotsc, z(t_k,r_k))$ as $n\to\infty$. 
The theorem is also valid in cases where one source of randomness has been 
turned off: if initial increments around $\tfl{n\ybar}$
 are deterministic then 
$v(\ybar)=0$, while if $D(\om)\equiv V$ then $\sigma_D^2=0$. 
The case $\sigma_D^2=0$ contains as special case the 
one with   deterministic weights $u^\om_\tau(k,j)\equiv p(0,j)$. 

If we consider only temporal correlations with a fixed
$r$, the formula for the  covariance is as follows: 
\be
\begin{split}
Ez(s,r)z(t,r)&= \frac{\kappa\sigma_a}{\sqrt{2\pi}}
 \vr(\ybar)^2\bigl( \sqrt{s+t\,}-\sqrt{ t-s}\;\bigr)\\
&\qquad\qquad
+ \frac{\sigma_a}{\sqrt{2\pi}}
v(\ybar)\bigl( \sqrt{s}+\sqrt{t}- \sqrt{s+t\,}\;\bigr)
\quad \text{ for $s<t$.} 
\end{split}
\label{rap-cov-2}
\ee

\begin{remark} 
The covariances are central to our proofs but they do not 
illuminate the behavior of the process $z$.  Here is a 
stochastic integral representation of the Gaussian process
with covariance \eqref{rap-cov}:
\be
\begin{split}
z(t,r)&=\vr(\ybar)\sigma_a\sqrt{\kappa}\iint_{[0,t]\times\Rb}
\varphi_{\sigma_a^2(t-s)}(r-x)\,dW(s,x)\\
&\qquad +\;  \sqrt{v(\ybar)}\int_{\Rb} 
\sign(x-r)
\Phi_{\sigma_a^2t}\bigl(-\abs{x-r}\,\bigr)\,dB(x).
\end{split}
\label{z-st-int}
\ee
Above $W$ is a two-parameter Brownian motion defined on 
$\Rb_+\times \Rb$,  $B$ is a one-parameter Brownian motion defined
on $\Rb$,  and $W$ and $B$ are independent of each other. 
The first integral represents the space-time noise created by the 
dynamics, and the second integral represents  the initial noise 
propagated by the evolution.  
The equality in \eqref{z-st-int} is equality  in distribution
of processes. It can  be verified by checking that the 
Gaussian process defined by the sum of the integrals has the covariance
\eqref{rap-cov}. 

One can readily see the second integral in \eqref{z-st-int} arise
as  a sum in the proof.  
It is the limit  of $Y^n(t,r)$ defined 
below equation \eqref{eq:hy}. 

One can also check that 
the right-hand side of \eqref{z-st-int}
 is a weak solution
of a stochastic heat equation with two independent sources of noise:
\be
z_t=\tfrac12{\sigma_a^2}\, z_{rr} +\vr(\ybar)\sigma_a\sqrt\kappa\, \dot{W}
+\tfrac12\sqrt{v(\ybar)}{\sigma_a^2}\,B'', 
\qquad z(0,r)\equiv 0.  
\label{z-st-eqn}
\ee
 $\dot{W}$ is space-time white noise generated by the dynamics
and $B''$  the second derivative
of the one-dimensional Brownian motion that represents initial noise. 
This equation has to be interpreted in a weak sense
through integration against smooth compactly supported
 test functions. We make a related remark  below in
Section \ref{rwre-lim-sec} for limit
processes of quenched means of space-time RWRE. 
\label{spde-remark1}
\end{remark}

The simplest RAP dynamics averages only two 
neighboring height values.  By translating the indices, we can 
assume that  $p(0,-1)+p(0,0)=1$. 
In this case the evolution of increments is given by
the equation 
\be
\eta_\tau(k)=u_\tau(k,0) \eta_{\tau-1}(k)+ u_\tau(k-1,-1)
\eta_{\tau-1}(k-1). 
\label{2pt-eta-evol}
\ee
There is a queueing interpretation of sorts for this 
evolution. Suppose $\eta_{\tau-1}(k)$ denotes the amount of work
that remains 
at station $k$ at the end of cycle $\tau-1$. 
Then during cycle $\tau$, 
the fraction $u_\tau(k,-1)$ of this work is completed 
and moves on to station $k+1$, while the remaining
fraction $u_\tau(k,0)$ stays at station $k$ for further processing. 

In this case  we can explicitly evaluate the 
constant $\beta$ in terms of the other quantities. 
In a particular stationary situation we can also 
identify the temporal marginal of $z$ in \eqref{rap-cov-2} 
 as a familiar process. (A probability distribution $\mu$ 
on the space  $\Zb^\Zb$ 
 is an invariant distribution for the increment process
if it is the case that when $\eta_0$ has $\mu$ distribution,
so  does $\eta_\tau$ 
for all times $\tau\in\Zb_+$.)

\begin{proposition} Assume  $p(0,-1)+p(0,0)=1$. 

{\rm (a)} Then 
\be
\beta= \frac1{\sigma_a^2}{\Ef [u_0(0,0)u_0(0,-1)]} .
\label{2pt-beta}
\ee

{\rm (b)} Suppose further that the increment process $\eta_\tau$
possesses
an invariant distribution $\mu$ in which the variables 
$\{\eta(i): i\in\Zb\}$ are i.i.d.~with common mean $\vr= 
E^\mu[\eta(i)]$ and variance 
$v=E^\mu[\eta(i)^2]-\vr^2$. Then $v=\kappa \vr^2$.

Suppose that in Theorem \ref{rap-thm-1} 
 each $\eta^n_\tau=\eta_\tau$ is a stationary process with marginal $\mu$. 
 Then  the limit process $z$ has  covariance 
\be
Ez(s,q)z(t,r) 
= \kappa\vr^2\Bigl( \Psi_{\sigma_a^2s} \bigl(\lvert q-r\rvert\bigr)+
\Psi_{\sigma_a^2t} \bigl(\lvert q-r\rvert\bigr)
-  \Psi_{\sigma_a^2\lvert t-s\rvert} \bigl(\lvert q-r\rvert\bigr)\Bigr).
\label{space-time-fBM}
\ee
In particular, for a fixed $r$ the process $\{z(t,r):t\in\Rb_+\}$
has covariance 
\be
Ez(s,r)z(t,r)=\frac{\sigma_a\kappa\vr^2}{\sqrt{2\pi}}
\bigl( \sqrt{s}+\sqrt{t}-\sqrt{\abs{t-s}}\,\bigr).
\label{fBM-1}
\ee
In other words, 
 process $z(\cdot,r)$ is fractional Brownian motion with Hurst parameter
$1/4$. 
\label{prop-2pt}
\end{proposition} 
 
To rephrase the connection \eqref{space-time-fBM}--\eqref{fBM-1},
 the process $\{z(t,r)\}$ in \eqref{space-time-fBM} is a certain 
two-parameter process whose marginals along the first parameter direction
are fractional Brownian motions.  

Ferrari and Fontes \cite{ferr-font-rap} showed that
given any slope $\rho$, the process $\eta_\tau$ started 
from deterministic increments $\eta_0(x)=\rho x$ 
converges weakly to an invariant distribution. 
But as is typical for interacting systems, there is 
little information about the invariant distributions in the general
case.  The next example  gives a family of processes and
i.i.d.~invariant distributions to show that part (b) of Proposition
\ref{prop-2pt} is not vacuous. Presently we are not aware of
other explictly known invariant distributions for RAP.

\begin{example} Fix integer parameters $m>j>0$. 
Let $\{u_\tau(k,-1): \tau\in\Nb, k\in\Zb\}$ 
be i.i.d.~ beta-distributed random variables
with density  
\[
h(u)=\frac{(m-1)!}{(j-1)!(m-j-1)!}u^{j-1}(1-u)^{m-j-1}
\]
on $(0,1)$. Set $u_\tau(k,0)=1-u_\tau(k,-1)$. Consider
the evolution defined by \eqref{2pt-eta-evol} with these weights.  
Then  a family of  invariant distributions for the
increment  process 
$\eta_\tau=(\eta_\tau(k): k\in\Zb)$ is obtained by
letting the variables $\{\eta(k)\}$ be i.i.d.~gamma 
distributed with common density 
\be
f(x)=\frac1{(m-1)!} \lambda e^{-\lambda x}(\lambda x)^{m-1}
\label{gamma-eq}
\ee
on $\Rb_+$. The family of invariant distributions is parametrized
by $0<\lambda<\infty$. Under this  distribution 
$
\Ev [\eta(k)]= {m}/\lambda 
$ and  
$\Vv[\eta(k)]={m}/{\lambda^2} . 
$
\label{example-1}
\end{example}

One motivation for the present work was to investigate
whether the limits found in \cite{sepp-rw} for 
fluctuations along a characteristic for 
independent walks are instances of some universal 
behavior. The present results are in agreement 
with those obtained for independent walks. 
The common scaling is $n^{1/4}$. In that paper only the 
case $r=0$ of Theorem \ref{rap-thm-1} was studied. 
For both independent walks and RAP the limit
$z(\cdot\,, 0)$  is a mean-zero Gaussian process 
with covariance of the type 
\[
Ez(s,0)z(t,0)= c_1\bigl( \sqrt{s+t\,}-\sqrt{ t-s}\;\bigr)
+ c_2\bigl( \sqrt{s}+\sqrt{t}- \sqrt{s+t\,}\;\bigr)
\]
where $c_1$ is determined by the mean increment 
and $c_2$ by the variance of the increment
locally  around
the initial point of the characteristic. Furthermore,
as in Proposition \ref{prop-2pt}(b), for independent  
walks the limit process specializes to fractional 
Brownian motion  if the increment process is stationary. 

These and other related results suggest several avenues
of inquiry. In the introduction we contrasted this
picture of $n^{1/4}$ fluctuations and fractional 
Brownian motion limits with the $n^{1/3}$ fluctuations
and Tracy-Widom limits found in exclusion and Hammersley
processes.    Obviously more classes of processes should
be investigated to understand better the demarcation
between these two types.  Also, there might be further
classes with different limits. 

Above we assumed independent  increments at time zero. It would be of
interest to see if relaxing this assumption leads to
a change in the second part of the covariance 
\eqref{rap-cov}. [The first part comes from the  
random walks in the dual description and would not be
affected by the initial conditions.]  However, without
knowledge of some explicit invariant distributions it is not clear
what types of initial increment processes  $\{\eta_0(k)\}$
are worth  considering.  Unfortunately finding 
explicit invariant  
distributions for interacting systems seems often  a
matter of good fortune.

We conclude this section with the dual description
of RAP which leads us to study random walks in
a space-time random environment. 
Given $\om$, let $\{X^{i,\tau}_s: s\in \Zb_+\}$ denote a random
walk on $\Zb$ that starts at $X^{i,\tau}_0=i$, and whose
transition probabilities are given by 
\be
P^\om(X_{s+1}^{i,\,\tau}=y\,|\,X^{i,\,\tau}_s=x)=u^\om_{\tau-s}(x,y-x).
\label{X-tr-pr}
\ee
$P^\om$ is the path measure of the walk $X^{i,\tau}_s$, with
expectation denoted by $E^\om$.  Comparison of 
 \eqref{si-dyn-1}  and \eqref{X-tr-pr} gives
 \be
\si_\tau(i)=\sum_j 
P^\om(X_{1}^{i,\,\tau}=j\,|\,X^{i,\,\tau}_0=i)  \si_{\tau-1}(j)
=E^\om\bigl[ \si_{\tau-1}(X_{1}^{i,\,\tau})\bigr].
\label{si-dyn-2}
\ee
Iteration and the Markov property of the walks $X^{i,\tau}_s$  then lead to 
\be
\si_\tau(i)=E^\om\bigl[\si_0(X^{i,\,\tau}_\tau)\bigr].
\label{eq:conn}
\ee
Note that the initial height function $\sigma_0$ is a constant 
under the expectation $E^\om$. 

Let us add another coordinate to keep track of time and write
$\Xbar^{i,\tau}_s=(X^{i,\tau}_s, \tau-s)$ for $s\geq 0$.  Then 
$\Xbar^{i,\tau}_s$ is a random walk on the planar lattice 
$\Zb^2$ that always moves down one step in the $e_2$-direction,
and if its current position is $(x,n)$,  the 
$e_1$-coordinate of its next position is $x+y$ with probability 
 $u_n(x, y)$.  We shall call it the   backward random
walk in a random environment.  In the next section we discuss
this walk and its   forward  counterpart. 
 
\section{Random walk in a space-time  random environment}
\label{rwre-section}
 
\subsection{Definition of the model}
We consider
here  a particular random walk in random environment (RWRE). 
The walk evolves on the planar integer lattice   $\Zb^2$, 
which we think 
of as space-time:
 the first component  
 represents one-dimensional discrete space,
and  the second   represents  discrete time. 
We denote by $e_2$ the unit vector in the time-direction. 
The walks will not be random in the $e_2$-direction, but only
in the spatial $e_1$-direction. 
 
 We consider {\sl forward walks} 
$\Zbar^{i,\tau}_m$ and 
{\sl backward walks} $\Xbar^{i,\tau}_m$.     The subscript $m\in\Zb_+$
is the time parameter of the walk
 and superscripts are initial points:
\be
\Zbar^{i,\tau}_0=\Xbar^{i,\tau}_0=(i,\tau)\in \Zb^2.
\label{init-pt}
\ee
The forward walks move deterministically up in time, while 
the backward walks move deterministically down in time:
\[
\Zbar^{i,\tau}_m=(Z^{i,\tau}_m, \tau+m)
\quad\text{and}\quad
\Xbar^{i,\tau}_m=(X^{i,\tau}_m, \tau-m)
\quad\text{for $m\geq 0$.} 
\] 
Since the time components of the walks are deterministic,
  only   the spatial components 
$Z^{i,\tau}_m$ and $X^{i,\tau}_m$ are really relevant. 
  We impose a finite 
range on the steps of the walks: there is a fixed constant $\rg$ such that 
\be
 \bigl\lvert Z^{i,\tau}_{m+1}- Z^{i,\tau}_m\bigr\rvert \leq \rg
\quad\text{and}\quad
 \bigl\lvert X^{i,\tau}_{m+1}- X^{i,\tau}_m\bigr\rvert \leq \rg. 
\label{rwre-bdrg}
\ee 

A note of advance justification for the setting:
The backward walks are the ones relevant to the random average process.
Distributions of forward and backward walks are obvious mappings
of each other.  However, we will be interested in the quenched mean
processes of the walks as we vary the final time for the 
forward walk or the initial space-time point for the backward walk. 
The results for the forward walk
 form an interesting point of comparison to the backward
walk, even though they will not be used to analyze the 
 random average process. 

An {\sl environment} 
  is a configuration  of   probability vectors 
$\om=\bigl(u_\tau(x): (x,\tau)\in\Zb^2\bigr)$  
where each vector 
$u_\tau(x)=(u_\tau(x,y): -\rg\leq y\leq \rg )$  satisfies 
 \[
\text{$0\leq u_\tau(x,y)\leq 1 $ for all $-\rg\leq y\leq \rg$, and }\; 
\sum_{y=-\rg}^{\rg}u_\tau(x,y)=1.  
\]
An environment $\om$ is a sample point  of the probability
space $(\Om,\Sk,\Pf)$. The sample space is the 
 product space
$\Om= \Pc^{\Zb^2}$ where $\Pc$ is the space of probability
vectors of length $2\rg+1$, and $\Sk$ is the product $\si$-field on $\Om$
induced by the Borel sets on $\Pc$.
Throughout, we assume that $\Pf$ is a product probability
 measure on $\Omega$
such that the vectors $\{ u_\tau(x):   (x,\tau)\in\Zb^2\}$ are
independent and identically 
distributed. 
  Expectation under $\Pf$ is denoted by $\Ef$.
When for notational convenience we wish to think of 
$u_\tau(x)$ as an infinite vector, then $u_\tau(x,y)=0$
for $\abs{y}>M$. We write $u^\om_\tau(x,y)$ to make explicit
the environment $\om$, and also $\om_{x,\tau}=u_\tau(x)$ for
the environment at space-time point $(x,\tau)$. 
 
Fix an environment  $\om$ and  an initial point $(i,\tau)$. 
The forward and backward  walks $\Zbar^{i,\tau}_m$ and 
  $\Xbar^{i,\tau}_m$  $(m\geq 0)$ 
 are defined as canonical $\Zb^2$-valued 
Markov chains  on their path spaces under the measure $P^\om$
determined by the conditions 
\[
\begin{split}
P^\om\{\Zbar^{i,\,\tau}_0=(i,\,\tau)\}&=1,\\
P^\om\{\Zbar_{s+1}^{i,\,\tau}=(y,\tau+s+1)\,|\,\Zbar^{i,\,\tau}_s=(x,\tau+s)\}
&=u_{\tau+s}(x,y-x) 
\end{split}
\]
for the forward walk, and by 
\[
\begin{split}
P^\om\{\Xbar^{i,\,\tau}_0=(i,\,\tau)\}&=1,\\
P^\om\{\Xbar_{s+1}^{i,\,\tau}=(y,\tau-s-1)\,|\,\Xbar^{i,\,\tau}_s=(x,\tau-s)\}
&=u_{\tau-s}(x,y-x) 
\end{split}
\]
for the backward walk. By dropping the time components
$\tau$, $\tau\pm s$ and $\tau\pm s\pm 1$ from the equations
 we get the corresponding 
properties for the spatial walks 
$Z^{i,\,\tau}_s$ and $X^{i,\,\tau}_s$. 
 When  we consider many walks 
under a common environment $\om$,   it will be notationally
convenient to attach the initial point $(i,\tau)$ to the walk
and only the environment $\om$ to the measure $P^\om$. 
 
  $P^\om$ is called the \emph{quenched} distribution, 
and expectation under $P^\om$ is denoted by $E^\om$. 
The \emph{annealed} distribution and  expectation 
are 
$
P(\cdot)=\Ef P^\om(\cdot)\quad\text{and}\quad E(\cdot)=\Ef E^\om(\cdot). 
$
  Under  $P$ both $X^{i,\tau}_m$  and $Z^{i,\tau}_m$ 
are ordinary  homogeneous random walks on $\Zb$  with jump probabilities 
$
p(i,\,i+j)=p(0,\,j)=\Ef u_0(0,j).   
$
  These walks satisfy the law of large numbers with velocity 
\be
V=\sum_{j\in\Zb}p(0,\,j)j.\label{eq:adrift}
\ee
As for RAP, we  also use the notation $b=-V$.  

\subsection{Limits for quenched mean processes}
\label{rwre-lim-sec}
We start by stating the 
 quenched invariance principle for   the space-time
RWRE.  $\{B(t):t\geq 0\}$ denotes 
standard one-dimensional Brownian motion. 
$D_{\Rb}[0,\infty)$ is the space of real-valued cadlag
functions on $[0,\infty)$ with the standard Skorohod metric
\cite{ethi-kurt}.  Recall the definition \eqref{def-si_a}
of the variance $\si_a^2$  of the annealed walk, and assumption
\eqref{ellipt} that guarantees that the quenched walk
has stochastic noise. 

\begin{theorem}\cite{qclt-spacetime} 
Assume {\rm \eqref{ellipt}}. Then for
  $\Pf$-almost every 
$\om$, under $P^\om$ the process 
$n^{-1/2}(X^{0,0}_{\tfl{nt}}-ntV)$ converges weakly    to
the process $B(\sigma_a^2t)$ on
the  path space  $D_\Rb[0,\infty)$ as $n\to\infty$. 

Assume further that $\si_D^2>0$ so that the environment 
is not degenerate.  Then
we have these bounds on the 
variance of the quenched mean: there exist 
constants $0<C_1,C_2<\infty$ such that for all $n$ 
\be
C_1n^{1/2}\leq \Ef\bigl[ (E^\om(X^{0,0}_n)-nV)^2\bigr]\leq C_2n^{1/2}. 
\label{qm-var-1}
\ee

\label{q-inv-pr}
\end{theorem}

Quite obviously, $X^{0,0}_n$ and $Z^{0,0}_n$ are 
interchangeable in the above theorem. 
Bounds \eqref{qm-var-1} suggest the possibility of a weak limit
for the quenched mean on the 
scale $n^{1/4}$.  Such results are the main point of this section.

For $t\geq 0$, $r\in\Rb$ we define scaled, centered quenched mean
processes 
\be
a_n(t,r)=n^{-1/4}\bigl\{ 
E^\om\bigl(Z^{\tfl{r\sqrt{n}}, 0}_{\tfl{nt}}\bigr)
 -\tfl{r\sqrt{n}}-\tfl{nt}V \bigr\} 
\label{def-an}
\ee
for the forward walks, and 
\be
y_n(t,r)=n^{-1/4}\bigl\{ 
E^\om\bigl(X^{\tfl{ntb}+\tfl{r\sqrt{n}}, \tfl{nt}}_{\tfl{nt}}\bigr)
 -\tfl{r\sqrt{n}} \bigr\}
\label{def-yn}
\ee
for the backward walks. In words, the process $a_n$ follows
forward  walks 
from level $0$ to level $\tfl{nt}$ and records
centered  quenched means. 
Process $y_n$ follows backward  walks 
from level $\tfl{nt}$ down to level $0$ and records the
centered  quenched mean 
of the point it hits at level 0. The initial points of the 
backward walks are translated by the negative of
the mean drift $\tfl{ntb}$. 
This way the temporal processes  $a_n(\cdot,r)$ and   $y_n(\cdot,r)$  obtained 
by fixing $r$ are meaningful processes. 

Random variable  $y_n(t,r)$ is not exactly
centered, for 
\be
\Ef y_n(t,r)=n^{-1/4}\bigl(\tfl{ntb}-\tfl{nt}b\bigr).
\label{yn-center}
\ee 
Of course this makes no difference to the limit.

Next we describe the Gaussian limiting processes. 
Recall the constant $\kappa$ defined in \eqref{def-kappa}
and the function $\Gamma_q$ defined
in \eqref{y-cov}. 
Let $\{a(t,r): (t,r)\in\Rb_+\times\Rb\}$ and
 $\{y(t,r): (t,r)\in\Rb_+\times\Rb\}$ be the mean zero Gaussian processes
with covariances
\[ Ea(s,q)a(t,r) =\Gamma_q\bigl((s\land t,q),(s\land t,r)\bigr) 
\]
and 
\[ Ey(s,q)y(t,r) =\Gamma_q\bigl((s,q),(t,r)\bigr) 
\] 
for $s, t\geq 0$ and $q,r\in\Rb$. 
When one argument is fixed, 
the random function  $r\mapsto y(t,r)$
is denoted by $y(t,\cdot)$ and $t\mapsto y(t,r)$
  by $y(\cdot, r)$.  From the covariances follows that 
at a fixed time level $t$ the spatial  processes  $a(t,\cdot)$ and 
$y(t,\cdot)$ are equal in distribution. 

We record basic properties of these processes. 

\begin{lemma}
The process $\{y(t,r)\}$ 
 has a version with continuous paths as functions
of $(t,r)$.  Furthermore, it has 
the following  Markovian structure
in time.  Given $0=t_0<t_1<\dotsm<t_n$, let 
$\{\ytil(t_i-t_{i-1},\cdot): 1\leq i\leq n\}$ be independent 
random functions such that   $\ytil(t_i-t_{i-1},\cdot)$ has
the distribution of $y(t_i-t_{i-1},\cdot)$ for $i=1,\dotsc,n$. 
Define $y^*(t_1,r)=\ytil(t_1,r)$ for $\,r\in\Rb$, and then inductively 
for $i=2,\dotsc,n$ and $r\in\Rb$,
\be
y^*(t_i,r)= \int_\Rb \varphi_{\sigma_a^2(t_i-t_{i-1})}(u)
y^*(t_{i-1},r+u) \,d u
+ \tilde{y}(t_i-t_{i-1}, r).
\label{y-markov-1}
\ee
Then the joint distribution of the random functions 
 $\{y^*(t_i,\cdot): 1\leq i\leq n\}$ is the same as that of
 $\{y(t_i,\cdot): 1\leq i\leq n\}$    from the original process. 
 \label{y-lemma-1}
\end{lemma}

\begin{proof}[Sketch of proof]
Consider $(s,q)$ and $(t,r)$ varying in a compact set. 
From the covariance comes the 
estimate
\be
E\bigl[(y(s,q)-y(t,r))^2\bigr] \leq C\bigl(\abs{s-t}^{1/2}+
 \abs{q-r}\bigr)
\label{y-mom-1}
\ee
from which, since the integrand is Gaussian, 
\be
E\bigl[(y(s,q)-y(t,r))^{10}\bigr] \leq C\bigl(\abs{s-t}^{1/2}+
 \abs{q-r}\bigr)^5 
\leq C\norm{(s,q)-(t,r)}^{5/2}.
\label{y-mom-2}
\ee
Kolmogorov's criterion implies the existence of a continuous
version. 

For the second  statement use \eqref{y-markov-1} to 
 express a linear 
combination $\sum_{i=1}^n \theta_i y^*(t_i,r_i)$ in
 the form 
\[
\sum_{i=1}^n \theta_i y^*(t_i,r_i)
= \sum_{i=1}^n 
\int_{\Rb} \ytil(t_i-t_{i-1}, x)\,\lambda_i(dx)
\]
where the signed measures  $\lambda_i$ are linear
combinations of Gaussian distributions.  Use this
representation to  compute
the variance of the linear combination on the left-hand
side (it is mean zero Gaussian).
Observe that this variance equals  
\[\sum_{i,j}\theta_i\theta_j\Gamma_q((t_i,r_i),(t_j,r_j)).
\qedhere \]
\end{proof}

\begin{lemma}
The process $\{a(t,r)\}$ 
 has a version with continuous paths as functions
of $(t,r)$.  Furthermore, it has 
independent increments 
in time. A more precise statement follows.
  Given $0=t_0<t_1<\dotsm<t_n$, let 
$\{\atil(t_i-t_{i-1},\cdot): 1\leq i\leq n\}$ be independent 
random functions such that   $\atil(t_i-t_{i-1},\cdot)$ has
the distribution of $a(t_i-t_{i-1},\cdot)$ for $i=1,\dotsc,n$. 
Define $a^*(t_1,r)=\atil(t_1,r)$ for $r\in\Rb$, and then inductively 
for $i=2,\dotsc,n$ and $r\in\Rb$,
\be
a^*(t_i,r)=  a^*(t_{i-1}  , r) + 
 \int_\Rb \varphi_{\sigma_a^2 t_{i-1}}(u)
\atil(t_i-t_{i-1} ,r+u) \,d u.
 \label{a-markov-1}
\ee
Then the joint distribution of the random functions 
 $\{a^*(t_i,\cdot): 1\leq i\leq n\}$ is the same as that of
 $\{a(t_i,\cdot): 1\leq i\leq n\}$    from the original process.
\label{a-lemma-1}
\end{lemma}

The proof of the lemma above is similar to the previous one
so we omit it.  

\begin{remark}  Processes $y$ and $a$ have 
representations in terms of stochastic  integrals.
As in Remark \ref{spde-remark1} let $W$ be a two-parameter Brownian motion on
$\Rb_+\times\Rb$. In more technical terms, $W$ is  the 
orthogonal Gaussian martingale  measure on $\Rb_+\times\Rb$
 with covariance 
$EW([0,s]\times A)W([0,t]\times B)=(s\land t)\Leb(A\cap B)$ 
for $s,t\in\Rb_+$
and  bounded Borel
sets $A,B\subseteq\Rb$. 
  Then 
\be
y(t,r)=\sigma_a\sqrt\kappa
\iint_{[0,t]\times\Rb} \varphi_{\sigma^2_a(t-s)}(r-z)\,
dW(s,z)
\label{y-st-int}
\ee
while
\be
a(t,r)=\sigma_a\sqrt\kappa\iint_{[0,t]\times\Rb} \varphi_{\sigma^2_as}(r-z)\,
dW(s,z).
\label{a-st-int}
\ee
By the equations above we mean equality in distribution of
 processes.  They  can be verified by a comparison of 
covariances, as the integrals on the right are also Gaussian
processes. 
Formula \eqref{y-st-int} implies
that  process $\{y(t,r)\}$ is a weak solution
of the stochastic heat equation
\be
y_t=\tfrac12{\sigma_a^2}\, y_{rr} +\sigma_a\sqrt\kappa\, \dot{W}\,,
\qquad y(0,r)\equiv 0 
\label{st-heat-eqn}
\ee
where $\dot{W}$ is white noise. 
(See \cite{wals}.)
These observations are not used elsewhere in the paper.
\end{remark}

Next we record the limits for the quenched mean processes. 
The four theorems that follow require assumption
\eqref{ellipt} of stochastic noise and the assumption that 
the annealed probabilities $p(0,j)=\Ef u_0^\om(0,j)$ have span 1.  
 This next theorem
is the one needed for Theorem \ref{rap-thm-1} for RAP.  

\begin{theorem}
The finite dimensional distributions of processes 
$y_n(t,r)$ converge to those of $y(t,r)$ as $n\to\infty$.
More precisely,
for any finite set of points $\{(t_j,r_j):1\leq j\leq k\}$ in
$\Rb_+\times\Rb$, 
the vector $\bigl(y_n(t_j,r_j):1\leq j\leq k\bigr)$ converges 
weakly in $\Rb^k$ to the vector 
$\bigl(y(t_j,r_j):1\leq j\leq k\bigr)$. 
\label{yn-thm-1}
\end{theorem}

Observe that 
property  \eqref{y-markov-1} is easy to understand from
the limit. It 
 reflects the Markovian property 
\[
E^\om(X^{x,\tau}_{\tau})=\sum_y P^\om(X^{x,\tau}_{\tau-s}=y)
E^\om(X^{y,s}_s) \quad \text{for $s<\tau$,} 
\]
  and the ``homogenization'' of the  coefficients which
converge to Gaussian probabilities by the quenched central 
limit theorem. 

Let us restrict the 
 backward quenched mean process to a single 
characteristic to observe the outcome. This is the source of 
the first term in the temporal correlations \eqref{rap-cov-2}
 for RAP.   The next statement  needs no proof, for it
  is just a particular case of the limit
in  Theorem \ref{yn-thm-1}.

\begin{corollary}
Fix $r\in\Rb$.  As $n\to\infty$, 
the finite dimensional distributions of the process 
$\{y_n(t,r): t\geq 0\}$ converge to those of the
mean zero Gaussian process
$\{y(t): t\geq 0\}$ with covariance
\[
Ey(s)y(t)
=\frac{\kappa\sigma_a}{\sqrt{2\pi}}
\bigl( \sqrt{t+s}-\sqrt{t-s}\,\bigr)
\qquad(s<t).
\]
\label{y-cor-1}
\end{corollary}

Then the same for the forward processes. 

\begin{theorem}   
The finite dimensional distributions of processes 
$a_n$ converge to those of $a$ as $n\to\infty$. More precisely,
for any finite set of points $\{(t_j,r_j):1\leq j\leq k\}$ in
$\Rb_+\times\Rb$, 
the vector $\bigl(a_n(t_j,r_j):1\leq j\leq k\bigr)$ converges 
weakly in $\Rb^k$ to the vector $\bigl(a(t_j,r_j):1\leq j\leq k\bigr)$. 
\label{an-thm-1}
\end{theorem}
 
When we specialize to a temporal 
process we also verify path-level tightness
and hence get weak convergence of the entire  process. 
When $r=q$ in (\ref{y-cov})
we get 
\[\Gamma_q\bigl((s\land t,r),(s\land t,r)\bigr) 
=c_a\sqrt{s\land t}\] with
 $c_a=
\sigma_D^2/(\beta \sqrt{\pi\sigma_a^2})$. 
Since $s\land t$ is the covariance of standard Brownian
motion $B(\cdot)$, we get the following limit.   

\begin{corollary}  Fix $r\in\Rb$.  
As $n\to\infty$, the   process
$\{a_n(t,r):t\geq 0\}$ converges weakly to 
$\{B(c_a\sqrt{t}\,):t\geq 0\}$  
 on the path space $D_{\Rb}[0,\infty)$. 
\label{a-cor-1}
\end{corollary}

\section{Random walk preliminaries} \label{prelim-sec} 
 In this section we collect some auxiliary results
for  random walks. 
The basic assumptions, \eqref{ellipt} and span 1
for the $p(0,j)=\Ef u_0(0,j)$ walk, 
 are in force
throughout the remainder of the paper. 

  Recall the  drift in the $e_1$ direction at the 
origin defined by 
\[
D(\om)=\sum_{x\in\Zb} x \; u^\om_0(0,x),
\]
with  
mean $V=-b=\Ef(D)$.   Define the centered drift by 
\[
g(\om)=D(\om)-V = E^\om(X^{0,0}_1-V). 
\]
The variance is   $\sigma_D^2=\Ef [g^2]$. 
The variance of the i.i.d.\ annealed walk in the $e_1$ direction 
 is 
\[
\sigma_a^2=\sum_{x\in\Zb} (x-V)^2\; \Ef u^\om_0(0,x).
\]
These variances are connected by
\[
\sigma_a^2 =\sigma_D^2 +  E[(X_1^{0,0}-D)^2].
\]
Let $X_n$ and $\tilde{X}_n$ be two independent walks in a
common environment $\om$, and $Y_n=X_n-\tilde{X}_n$. 
In the annealed sense
$Y_n$ is a Markov chain on $\Zb$ with transition probabilities 
\begin{align*}
q(0,y)&=\sum_{z\in\Zb} \Ef[u_0(0,z)u_0(0,z+y)] \qquad (y\in\Zb)\\
q(x,y)&=\sum_{z\in\Zb} p(0,z)p(0,z+y-x)
\qquad (x\neq 0,y\in\Zb).
\end{align*}
$Y_n$ can be thought of 
as  a symmetric random walk on $\Zb$ whose transition
has been perturbed at the origin. The corresponding  homogeneous,
unperturbed 
transition probabilities are 
\[
\bar{q}(x,y)=\bar{q}(0,y-x)=
\sum_{z\in\Zb} p(0,z)p(0,z+y-x) 
\qquad (x ,y\in\Zb).
\]
The $\qbar$-walk has variance $2\sigma_a^2$ and 
 span 1 as can be deduced from the definition and 
the hypothesis that the $p$-walk has span 1.
 Since the $\qbar$-walk is symmetric, its range must be a
subgroup of $\Zb$. Then span 1 implies  that 
it is irreducible. The $\qbar$-walk is recurrent by the Chung-Fuchs theorem.
Elementary arguments extend 
  irreducibility and recurrence from $\qbar$ 
  to the $q$-chain because away from the origin the two
walks are the same. 
Note that assumption \eqref{ellipt} is required here because
  the $q$-walk is absorbed at the origin iff \eqref{ellipt} fails. 
 
Note that the functions defined in \eqref{def-chf-2} are 
the characteristic functions of these transitions:
\[
\lambda(t)= \sum_x q(0,x)e^{itx}
\quad\text{and}\quad  
\bar{\lambda}(t)= \sum_x \qbar(0,x)e^{itx}.
\]
 
Multistep transitions are denoted by  $q^k(x,y)$ and 
$\bar{q}^k(x,y)$, defined as usual  by  
\begin{align*}
&\qbar^0(x,y)=\ind_{\{x=y\}},\quad \bar{q}^1(x,y)=q(x,y),  \\
&\qquad \qquad \qbar^k(x,y)=\sum_{x_1,\dotsc,x_{k-1}\in\Zb}
 \qbar(x,x_1)\qbar(x_1,x_2)\dotsm \qbar(x_{k-1},y) \quad (k\geq 2).
\end{align*}
 Green functions for the $\qbar$- and $q$-walks  are 
\[
\bar{G}_n(x,y)=\sum_{k=0}^n \bar{q}^k(x,y)
\quad\text{and}\quad
{G}_n(x,y)=\sum_{k=0}^n {q}^k(x,y).
\]
 $\bar{G}_n$ is symmetric but $G_n$ not necessarily. 

The potential kernel $\abar$ of the $\qbar$-walk is
defined by 
\be
\abar(x)=\lim_{n\to\infty} \bigl\{ \Gbar_n(0,0)-\Gbar_n(x,0)\bigr\}.
\label{abar-def}
\ee
It satisfies $\abar(0)=0$, 
  the equations 
\be
\abar(x)=\sum_{y\in\Zb} \qbar(x,y)\abar(y) 
\quad\text{for $x\neq 0$, and} \quad 
\sum_{y\in\Zb} \qbar(0,y)\abar(y)=1, 
\label{abar-eq}
\ee
and the limit 
\be 
\lim_{x\to\pm\infty}\frac{\abar(x)}{\abs{x}} = \frac1{2\sigma_a^2}.
\label{abar-lim}
\ee
These facts can be found in Sections 28 and 29 of Spitzer's
monograph \cite{spit}.  

\begin{example} If for some $k\in\Zb$, $p(0,k)+p(0,k+1)=1$, so that
$\qbar(0,x)=0$ for $x\notin\{-1,0,1\}$, then 
$\abar(x)=\abs{x}/(2\sigma_a^2)$. 
\label{ex-abar-1}
\end{example}

Define the constant 
\be
\beta=\sum_{x\in\Zb} q(0,x)\abar(x).
\label{def-beta-2}
\ee
To see that this definition agrees with \eqref{def-beta-1}, 
observe that the above equality leads to 
\begin{align*}
\beta 
= \lim_{n\to\infty} \Bigl\{\; \sum_{k=0}^n  \qbar^k(0,0)
- \sum_{k=0}^n \sum_x q(0,x)\qbar^k(x,0)  \;\Bigr\}.
\end{align*}
Think of the last sum over $x$ as $P[Y_1+\Ybar_k=0]$ 
where $\Ybar_k$ is the $\qbar$-walk, and 
$Y_1$ and $\Ybar_k$ are
independent. Since 
$Y_1+\Ybar_k$ has characteristic function 
$\lambda(t)\bar{\lambda}^k(t)$, we get 
\begin{align*}
\beta=
\lim_{n\to\infty}   \frac1{2\pi} \int_{-\pi}^\pi
(1-\lambda(t))\sum_{k=0}^n \bar{\lambda}^k(t)\,dt 
= \frac1{2\pi} \int_{-\pi}^\pi
 \frac{1- \lambda(t)} {1- \bar{\lambda}(t)} \,dt.
\end{align*}

Ferrari and Fontes \cite{ferr-font-rap} 
begin their development by showing that 
\[
\beta =\lim_{s\nearrow 1} \frac{\bar{\zeta}(s)}{\zeta(s)}
\]
where $\zeta$ and $\bar{\zeta}$ are the  generating functions 
\[
\zeta(s) =\sum_{k=0}^\infty q^k(0,0)s^k
\quad\text{and}\quad
\bar{\zeta}(s) =\sum_{k=0}^\infty \bar{q}^k(0,0)s^k.
\]

Our development bypasses  the generating functions. 
We begin with the asymptotics of the Green functions. 
This is the key to all our results, both for RWRE and RAP. 
As already pointed out, without assumption \eqref{ellipt}
the result would be completely wrong because the 
$q$-walk absorbs at 0, while a span $h>1$ would appear
in this limit as an extra factor. 
 
\begin{lemma} Let $x\in\Rb$, and 
let $x_n$ be any sequence of integers such that 
$x_n- n^{1/2}x$ stays bounded. 
Then
\be
\lim_{n\to\infty} 
{n}^{-1/2} {G}_n\bigl(x_n,0\bigr)
=\frac{1}{2\beta\sigma_a^2} \int_0^{2\sigma_a^2} 
\frac1{\sqrt{2\pi v}}
\exp\Bigl\{-\frac{x^2}{2v}\Bigr\}  \,d v. 
\label{G-lim1}
\ee
\end{lemma}

\begin{proof}
For the homogeneous $\qbar$-walk
the local limit theorem \cite[Section 2.5]{durr} 
implies that 
\be
\lim_{n\to\infty} 
{n}^{-1/2} \bar{G}_n(0,x_n)
=\frac1{2\sigma_a^2} \int_0^{2\sigma_a^2} \frac1{\sqrt{2\pi v}}
\exp\Bigl\{-\frac{x^2}{2v}\Bigr\}  \,d v
\label{barG-lim1}
\ee
and by symmetry the same limit is true for 
${n}^{-1/2} \bar{G}_n(x_n,0)$.
In particular, 
\be
\lim_{n\to\infty} 
{n}^{-1/2} \bar{G}_n(0,0)
=\frac1{\sqrt{\pi\sigma_a^2}}\,.
\label{barG-lim2}
\ee

Next we show 
\be
\lim_{n\to\infty} 
{n}^{-1/2} {G}_n(0,0)
=\frac1{\beta\sqrt{ \pi\sigma_a^2}}.  
\label{barG-lim2.5}
\ee
Using \eqref{abar-eq}, $\abar(0)=0$,
 and $\qbar(x,y)=q(x,y)$ for $x\neq 0$ we 
develop
\begin{align*}
&\sum_{x\in\Zb} q^m(0,x)\abar(x) =\sum_{x\neq 0} q^m(0,x)\abar(x) =
\sum_{x\neq 0, y\in\Zb}
 q^m(0,x) \qbar(x,y)\abar(y)\\
& =\sum_{x\neq 0, y\in\Zb}
 q^m(0,x) q(x,y)\abar(y)
=\sum_{y\in\Zb} q^{m+1}(0,y)\abar(y) 
-q^m(0,0) \sum_{y\in\Zb} q(0,y)\abar(y).
\end{align*}
Identify $\beta$ in the last sum above and 
  sum over $m=0,1,\dotsc, n-1$ to get 
\[
\bigl(1+q(0,0)+\dotsm +q^{n-1}(0,0)\bigr)\beta
= \sum_{x\in\Zb} q^n(0,x)\abar(x).
\]
Write this in the form 
\[
{n}^{-1/2} {G}_{n-1}(0,0)\beta = {n}^{-1/2} E_0\bigl[\abar(Y_n)\bigr].
\]
Recall that $Y_n=X_n-\Xtil_n$ where $X_n$ and $\Xtil_n$ are two
independent walks in the same environment.  Thus by Theorem
\ref{q-inv-pr} ${n}^{-1/2}Y_n$ converges weakly to 
a centered Gaussian with variance $2\sigma_a^2$. Under the 
annealed measure the walks $X_n$ and $\Xtil_n$ are 
ordinary i.i.d.\ walks with bounded steps, hence there is 
enough uniform integrability to conclude that
${n}^{-1/2} E_0\abs{Y_n} \to 2\sqrt{\sigma_a^2/\pi}.$
By \eqref{abar-lim} and straightforward estimation, 
\[
{n}^{-1/2} E_0\bigl[\abar(Y_n)\bigr] \to \frac1{\sqrt{\sigma_a^2\pi}}\,.
\]
This proves \eqref{barG-lim2.5}.

From \eqref{barG-lim2}--\eqref{barG-lim2.5}  we take the conclusion 
\be
\lim_{n\to\infty} \frac1{\sqrt{n}} 
\bigl\lvert {\beta} {G}_n(0,0)
\,-\,   \bar{G}_n(0,0) \bigr\rvert =0.
\label{barG-lim3}
\ee

Let $f^0(z,0)=\ind_{\{z=0\}}$ and for $k\geq 1$ let
\[
f^k(z,0)=\ind_{\{z\neq 0\}} \sum_{z_1\neq 0,\dotsc,z_{k-1}\neq 0} 
q(z,z_1)q(z_1,z_2)\dotsm q(z_{k-1},0).
\]
This is the probability that the first visit to the origin
 occurs at time $k$, including a possible first
visit at time $0$.  Note that 
this quantity is the same for the $q$ and $\bar{q}$ walks. 
Now bound 
\begin{align*}
&\sup_{z\in\Zb}\; \Bigl\lvert \frac{\beta}{\sqrt{n}} {G}_n(z,0)
\,-\, \frac1{\sqrt{n}} \bar{G}_n(z,0) \Bigr\rvert \\
&\quad \leq \sup_{z\in\Zb}\; \frac1{\sqrt{n}} \sum_{k=0}^n f^k(z,0) 
\bigl\lvert \beta  {G}_{n-k}(0,0)
\,-\,  \bar{G}_{n-k}(0,0) \bigr\rvert.
\end{align*}
To see that the last line vanishes as $n\to\infty$, 
by \eqref{barG-lim3} choose $n_0$ so that 
\[\lvert \beta  {G}_{n-k}(0,0)\,-\,
 \bar{G}_{n-k}(0,0) \rvert\leq \e\sqrt{n-k}\] for $k\leq n-n_0$,
while trivially 
\[\lvert \beta  {G}_{n-k}(0,0)\,-\,
 \bar{G}_{n-k}(0,0) \rvert \leq Cn_0\] for $n-n_0<k\leq n$.
The conclusion \eqref{G-lim1} now follows from this and  \eqref{barG-lim1}.
 \end{proof}

\begin{lemma} 
$\displaystyle{\sup_{n\geq 1} \,\sup_{x\in\Zb} \,
\bigl\lvert G_n(x,0)-G_n(x+1,0)\bigr\rvert
<\infty}.$
\label{Gn-lm-2}
\end{lemma}

\begin{proof}
Let $T_{y}=\inf\{n\geq 1: Y_n=y\}$ 
denote the first hitting time of the point $y$. 
\begin{align*}
G_n(x,0)&=E_x\Bigl[\;\sum_{k=0}^n \ind\{Y_k=0\}\Bigr]
=E_x\Bigl[\;\sum_{k=0}^{T_y\land n} \ind\{Y_k=0\}\Bigr]
+E_x\Bigl[\;\sum_{k=T_y\land n+1}^n \ind\{Y_k=0\}\Bigr]\\
&\leq E_x\Bigl[\;\sum_{k=0}^{T_y} \ind\{Y_k=0\}\Bigr]
+G_n(y,0).
\end{align*} 
In an irreducible  Markov chain the expectation 
$E_x\bigl[\;\sum_{k=0}^{T_y} \ind\{Y_k=0\}\bigr]$ is finite 
for any given states $x,y$ \cite[Theorem 3 in Section I.9]{chung}. 
 Since this is independent of $n$,
the inequalities above show that 
\be
\sup_n\sup_{-a\leq x\leq a}\lvert G_n(x,0)-G_n(x+1,0)\rvert
<\infty
\label{Gn-temp-2}
\ee
for any fixed $a$.

Fix a positive integer $a$ larger than the range of the jump
kernels $q(x,y)$ and $\qbar(x,y)$. Consider $x>a$. Let
$\sigma=\inf\{n\geq 1: Y_n\leq a-1\}$ 
and $\tau=\inf\{n\geq 1: Y_n\leq a\}$. Since the $q$-walks starting 
at $x$ and $x+1$ obey the translation-invariant kernel 
$\qbar$ until they hit the origin, 
\[
P_x[Y_{\sigma}=y, \sigma= n] =P_{x+1}[Y_{\tau}=y+1,\tau= n].  
\]
(Any path that starts at $x$ and enters $[0,a-1]$ 
at $y$ can be translated by 1 to a path that starts at $x+1$
and enters $[0,a]$ at $y+1$, without changing its probability.) 
Consequently 
\begin{align*}
 G_n(x,0)-G_n(x+1,0) 
= \sum_{k=1}^n \sum_{y=0}^{a-1} P_x[Y_{\sigma}=y, \sigma=k]
\bigl( G_{n-k}(y,0)-G_{n-k}(y+1,0)\bigr). 
\end{align*}
Together with \eqref{Gn-temp-2}
this shows that the quantity  in the statement of the lemma
is uniformly bounded over $x\geq 0$. The same argument works for $x\leq 0$. 
\end{proof}

One can also derive the limit
\[
\lim_{n\to\infty} \bigl\{  G_n(0,0)-G_n(x,0)\bigr\}= \beta^{-1}\abar(x)
\]
but we have no need for this. 

Lastly,  a moderate deviation bound 
for the  space-time RWRE with bounded steps.
Let  $X^{i,\tau}_s$ be   the spatial backward
walk    defined in 
Section \ref{rwre-section} with the bound \eqref{rwre-bdrg} on the steps. 
Let  $\wt X^{i,\,\tau}_s=X^{i,\,\tau}_s-i-Vs$ be the centered walk. 

\begin{lemma}\label{lm:moddev} 
For $m,n\in\Nb$, let $(i(m,n),\tau(m,n))\in\Zb^2$, $v(n)\geq 1$,
 and let
 $s(n)\to\infty$  be a sequence 
of positive integers.  Let $\al$, $\ga$  and $c$ be positive reals. 
Assume 
\[
\sum_{n=1}^\infty v(n)s(n)^\al \exp\{-cs(n)^{\ga}\}<\infty.
\]
Then for $\Pf$-almost every $\om$, 
\begin{equation}
\lim_{n\to\infty} \;\max_{1\leq m\leq v(n)}\;
s(n)^\al
 P^\om\Bigl\{\max_{1\leq k\leq s(n)}
\wt X^{i(m,n),\,\tau(m,n)}_k\geq cs(n)^{\frac12+\ga}
\Bigr\}=0.
\end{equation}
 \end{lemma}
\begin{proof}
Fix $\ve>0$. By Markov's inequality  and translation-invariance, 
\[
\begin{split}
&\Pf\biggl[\om\,:\;  \max_{1\leq m\leq v(n) }\; s(n)^\al
 P^\om\Bigl\{
\max_{1\leq k\leq s(n)}
\wt X^{i(m,n),\,\tau(m,n)}_k\geq cs(n)^{\frac12+\ga}
\Bigr\}  \geq \ve\biggr] \\
&\quad \leq  {\ve}^{-1}{s(n)^\al} v(n)
P\bigl\{\max_{1\leq k\leq s(n)}
\wt X^{0,0}_k\geq cs(n)^{\frac12+\ga}
\bigr\}. 
 \end{split}
\]
 Under the annealed measure $P$,   $\wt X^{0,0}_k$ is an ordinary homogeneous
mean zero  random walk with bounded steps. 
It has a finite moment generating function
$\phi(\lambda)=\log E(\exp\{\lambda \wt X^{0,0}_1\})$ that satisfies
$\phi(\lambda)=\Oc(\lambda^2)$ for small $\lambda$. 
Apply  Doob's inequality to the martingale 
$M_k=\exp(\lambda \wt X^{0,0}_k -k\phi(\lambda))$, note that
$\phi(\lambda)\geq 0$, and choose  a  constant $a_1$  such that
 $\phi(\lambda)\leq a_1\lambda^2$ for small $\lambda$.
 This gives 
\begin{align*}
&P\bigl\{\max_{1\leq k\leq s(n)}
\wt X^{0,0}_k\geq cs(n)^{\frac12+\ga}
\bigr\}
\leq 
P\bigl\{\max_{1\leq k\leq s(n)} M_k
\geq \exp\bigl(c\lambda s(n)^{\frac12+\ga}-s(n)\phi(\lambda)\bigr) \bigr\}\\
&\quad \leq \exp\bigl(-c\lambda s(n)^{\frac12+\ga}+a_1s(n)\lambda^2\bigr)
= e^{a_1}\cdot \exp\{-cs(n)^{\ga}\}
 \end{align*}
where we took $\lambda=  s(n)^{-\frac12}$.

The conclusion of the lemma now 
 follows from the hypothesis and Borel-Cantelli. 
\end{proof}

\section{Proofs for backward walks in a random environment}
\label{bw-walk-pf-sec} 
Here are two further notational conventions used in the proofs. 
The environment configuration at a fixed time level 
is denoted by $\ombar_n=\{\om_{x,n}: x\in\Zb\}$. 
 Translations on $\Omega$
are defined by $(T_{x,n}\om)_{y,k}=\om_{x+y,n+k}$.

\subsection{Proof of Theorem \ref{yn-thm-1}}
This proof  proceeds in two stages.
First in Lemma \ref{yn-lemma-1} convergence is proved 
for finite-dimensional distributions at a fixed $t$-level.
In the second stage the convergence is extended to 
multiple $t$-levels via the natural Markovian property
that we express in terms of $y_n$  next.
Abbreviate 
$X^{n,t,r}_k=
X^{\tfl{ntb}+\tfl{r\sqrt{n}}, \tfl{nt}}_k$. Then for
$0\leq s<t$,
\begin{align}
&y_n(t,r)=  n^{-1/4}\bigl(
E^\om(X^{n,t,r}_{\tfl{nt}})-\tfl{r\sqrt{n}}\bigr) \nn\\
&=\sum_{z\in\Zb} P^\om\bigl\{ X^{n,t,r}_{\tfl{nt}-\tfl{ns}}
=\tfl{nsb}+  z\bigr\}
 n^{-1/4}\bigl(
E^\om(X^{\tfl{nsb}+  z, \tfl{ns}}_{\tfl{ns}})
-z\bigr)\nn\\
&\qquad  + \sum_{z\in\Zb} P^\om\bigl\{ X^{n,t,r}_{\tfl{nt}-\tfl{ns}}
=\tfl{nsb}+  z\bigr\}
 n^{-1/4}\bigl( z -\tfl{r\sqrt{n}}\bigr)   \nn\\
&=\sum_{z\in\Zb} P^\om\bigl\{ X^{n,t,r}_{\tfl{nt}-\tfl{ns}}
=\tfl{nsb}+  z\bigr\}
 n^{-1/4}\bigl(
E^\om(X^{\tfl{nsb}+  z, \tfl{ns}}_{\tfl{ns}})
-z\bigr)\nn\\
&\qquad + 
n^{-1/4}\bigl\{ E^\om\bigl( X^{n,t,r}_{\tfl{nt}-\tfl{ns}}\bigr)
-\tfl{nsb} -\tfl{r\sqrt{n}}\bigr\}  \nn\\
&=\sum_{z\in\Zb} P^\om\bigl\{ X^{n,t,r}_{\tfl{nt}-\tfl{ns}}
=\tfl{nsb}+  z\bigr\}
 n^{-1/4}\bigl(
E^\om(X^{\tfl{nsb}+  z, \tfl{ns}}_{\tfl{ns}})
-z\bigr)\label{yn-decomp-l-1}\\
&\qquad + 
y_n(u_n,r)\circ T_{\tfl{ntb}-\tfl{nbu_n}, \tfl{nt}-\tfl{nu_n}}
     + n^{-1/4}\bigl( \tfl{ntb}-\tfl{nsb}-\tfl{nbu_n} \bigr)
\label{yn-decomp-l-2}
\end{align}
where we defined $u_n=n^{-1}(\tfl{nt}-\tfl{ns})$ so that 
$\tfl{nu_n}=\tfl{nt}-\tfl{ns}$. 
$ T_{x,m}$ denotes the translation of the random environment
that makes  $(x,m)$ the new space-time origin, in other words  
$( T_{x,m}\om)_{y,n}=\om_{x+y, m+n}$.
 
The key to making use of the decomposition of $y_n(t,r)$ given 
on lines \eqref{yn-decomp-l-1} and \eqref{yn-decomp-l-2} is that 
the quenched expectations 
\[E^\om(X^{\tfl{nsb}+  z, \tfl{ns}}_{\tfl{ns}})\quad
\text{and}\quad 
y_n(u_n,r)\circ T_{\tfl{ntb}-\tfl{nbu_n}, \tfl{nt}-\tfl{nu_n}}\]
are independent because they are functions of 
environments $\ombar_m$ on disjoint sets of levels $m$, while 
the coefficients $P^\om\bigl\{ X^{n,t,r}_{\tfl{nt}-\tfl{ns}}
=\tfl{nsb}+  z\bigr\}$ on line \eqref{yn-decomp-l-1}  
converge (in probability) to Gaussian probabilities by 
the quenched CLT as $n\to\infty$.  In the limit this
decomposition becomes 
\eqref{y-markov-1}. 

Because of the little technicality of matching 
$\tfl{nt}-\tfl{ns}$ with $\tfl{n(t-s)}$ we state the 
next lemma for a sequence $t_n\to t$ instead of a fixed $t$.

\begin{lemma} Fix $t>0$, and finitely many reals 
$r_1<r_2<\dotsc<r_N$.  Let $t_n$ be a sequence
of positive reals such that  $t_n\to t$. 
 Then as $n\to\infty$  the $\Rb^N$-valued vector 
$(y_n(t_n,r_1),\dotsc, y_n(t_n,r_N))$ converges weakly to a 
mean zero Gaussian vector with covariance matrix 
$\{\Gamma_q((t,r_i),(t,r_j)): 1\leq i,j\leq N\}$ with 
$\Gamma_q$ as defined in \eqref{y-cov}.
\label{yn-lemma-1}
\end{lemma}

The proof of Lemma \ref{yn-lemma-1} is technical 
(martingale CLT and random walk estimates), so we postpone
it and proceed with the main development. 

\begin{proof}[Proof of Theorem \ref{yn-thm-1}]
The argument is inductive on the number $M$ of time points 
in the finite-dimensional distribution. 
 The induction assumption is that  
\be
\begin{split}
&\text{$[y_n(t_i,r_j): 1\leq i\leq M, 1\leq j\leq N] \to 
[y(t_i,r_j): 1\leq i\leq M, 1\leq j\leq N]$}  \\
&\text{weakly on $\Rb^{MN}$ for any $M$ time points 
$0\leq t_1<t_2<\dotsm<t_{M}$ and }\\
&\text{for any reals  $r_1,\dotsc,r_N$ for any finite $N$.}
\end{split}
\label{yn-ind-M}
\ee
The case $M=1$ comes from Lemma \ref{yn-lemma-1}. 
To handle the case $M+1$, let
$0\leq t_1<t_2<\dotsm<t_{M+1}$, and 
fix an arbitrary $(M+1)N$-vector 
$[ \theta_{i,j}]$. By the Cram\'er-Wold device, 
it suffices to show the weak convergence of the 
linear combination
\be
\sum_{\substack{1\leq i\leq M+1\\[2pt] 1\leq j\leq N}} \theta_{i,j}
y_n(t_i,r_j) 
=\sum_{\substack{1\leq i\leq M\\[2pt] 1\leq j\leq N}} \theta_{i,j}
y_n(t_i,r_j) 
+ \sum_{ 1\leq j\leq N} \theta_{M+1,j}
y_n(t_{M+1},r_j) 
\label{yn-M-decomp}
\ee
where we separated out the $(M+1)$-term to be manipulated.
The argument will use \eqref{yn-decomp-l-1}--\eqref{yn-decomp-l-2}
to replace the values at $t_{M+1}$ with values at $t_{M}$ plus
terms independent of the rest. 

For Borel sets  $B\subseteq\Rb$ define the probability measure 
\[
p^\om_{n,j}(B)= 
P^\om\bigl\{ X^{\tfl{nt_{M+1}b}+\tfl{r_j\sqrt{n}}, 
\tfl{nt_{M+1}}}_{\tfl{nt_{M+1}}-\tfl{nt_{M}}}-\tfl{nt_Mb}
\in B\bigr\}.
\]
Apply the decomposition 
\eqref{yn-decomp-l-1}--\eqref{yn-decomp-l-2}, with 
 $s_n=n^{-1}(\tfl{nt_{M+1}}-\tfl{nt_M})$ and 
\[
\tilde{y}_n(s_n,r_j) = 
y_n(s_n,r_j)\circ T_{\tfl{nt_{M+1}b}-\tfl{ns_nb}, 
\tfl{nt_{M+1}}-\tfl{ns_n}}
\]
to get 
\be
\begin{split}
y_n(t_{M+1},r_j)
&=\sum_{z\in\Zb} p^\om_{n,j}(z)  n^{-1/4}\bigl\{
E^\om(X^{\tfl{nt_Mb}+  z, \tfl{nt_M}}_{\tfl{nt_M}})
-z\bigr\} \\
&\qquad\qquad\qquad\qquad  + \tilde{y}_n(s_n,r_j) +\Oc(n^{-1/4}).
\end{split} 
\label{yn-temp-4}
\ee
The $\Oc(n^{-1/4})$ term above is 
$n^{-1/4}\bigl( \tfl{nt_{M+1}b}-\tfl{nt_Mb}-\tfl{ns_nb} \bigr)$,
 a deterministic quantity. 
Next we reorganize the sum in \eqref{yn-temp-4} to take
advantage of Lemma \ref{yn-lemma-1}. 
Given $a>0$, define a partition of $[-a,a]$ by 
\[-a=u_0<u_1<\dotsm<u_L=a\]    with mesh 
$\Delta=\max\{u_{\ell+1}-u_\ell\}$. 
For integers $z$ such that  $-a\sqrt{n}<z\leq a\sqrt{n}$,
let $u(z)$ denote the value $u_\ell$ such that  
$u_\ell\sqrt{n}<z\leq u_{\ell+1}\sqrt{n}$. 
For $1\leq j\leq N$ define an error term by 
\begin{align}
R_{n,j}(a) &=
n^{-1/4}\sum_{z=\tfl{-a\sqrt{n}}+1}^{\tfl{a\sqrt{n}}}
p^\om_{n,j}(z)
\Bigl( \bigl\{
E^\om(X^{\tfl{nt_Mb} +  z, \tfl{nt_M}}_{\tfl{nt_M}})
-z\bigr\}  \nn\\
& \qquad\qquad -\;  \bigl\{
E^\om(X^{\tfl{nt_Mb}+  \tfl{u(z)\sqrt{n}},
 \tfl{nt_M}}_{\tfl{nt_M}})
-\tfl{u(z)\sqrt{n}}\bigr\} \,\Bigr)
\label{def-R(a)-1}\\
&+ n^{-1/4}\sum_{z\leq -a\sqrt{n}\,,\; z> a\sqrt{n}}
p^\om_{n,j}(z)
 \bigl\{
E^\om(X^{\tfl{nt_Mb} +  z, \tfl{nt_M}}_{\tfl{nt_M}})
-z\bigr\}.
\label{def-R(a)-2}
\end{align}
With this we can rewrite  \eqref{yn-temp-4} as 
\be
\begin{split}
y_n(t_{M+1},r_j)
&=\sum_{\ell=0}^{L-1} p^\om_{n,j}(u_\ell{n}^{1/2}, 
u_{\ell+1}{n}^{1/2}\,]
  y_n(t_M, u_\ell) + \tilde{y}_n(s_n,r_j)  \\
&\qquad\qquad\qquad\qquad +R_{n,j}(a)  +\Oc(n^{-1/4}).
\end{split} 
\label{yn-temp-5}
\ee
Let $\gamma$ denote a normal distribution on $\Rb$ with mean
zero and variance $\sigma_a^2(t_{M+1}-t_M)$. 
According to the quenched CLT Theorem \ref{q-inv-pr}, 
\be
p^\om_{n,j}(u_\ell{n}^{1/2}, 
u_{\ell+1}{n}^{1/2}\,]\to \gamma(u_\ell-r_j, u_{\ell+1}-r_j]
\quad\text{in $\Pf$-probability  as $n\to\infty$. }
\label{q-clt-2}
\ee

In view of \eqref{yn-M-decomp} and \eqref{yn-temp-5}, 
we can write 
\be
\begin{split} 
\sum_{\substack{1\leq i\leq M+1\\[2pt] 1\leq j\leq N}} \theta_{i,j}
y_n(t_i,r_j) 
&=\sum_{\substack{1\leq i\leq M\\[2pt] 1\leq k\leq K}} \rho^\om_{n,i,k}\,
y_n(t_i,v_k) 
+
\sum_{1\leq j\leq N} \theta_{M+1,j}\,
\tilde{y}_n(s_n,r_j)\\
&\qquad   + R_n(a) +\Oc(n^{-1/4}). 
\end{split}
\label{yn-temp-6}
\ee
Above the spatial points $\{v_k\}$ are a relabeling of
$\{r_j, u_\ell\}$, the  $\om$-dependent coefficients 
$\rho^\om_{n,i,k}$ contain constants  $\theta_{i,j}$, probabilities
$p^\om_{n,j} (u_\ell{n}^{1/2}, 
u_{\ell+1}{n}^{1/2}\,]$, and zeroes. 
The constant limits $\rho^\om_{n,i,k}\to\rho_{i,k}$ exist
in $\Pf$-probability as $n\to\infty$. 
 The error in \eqref{yn-temp-6} is 
 $R_n(a)=\sum_j \theta_{M+1,j} R_{n,j}(a)$. 

The variables $\tilde{y}_n(s_n,r_j)$ are functions of
the environments 
$\{\ombar_m: [nt_{M+1}]\geq m> [nt_M]\}$ and hence
independent of $y_n(t_i,v_k)$ for $1\leq i\leq M$ which 
 are functions of
$\{\ombar_m: [nt_M]\geq m> 0\}$.   

On a probability space on which the limit process $\{y(t,r)\}$ has been
defined, let 
 $\tilde{y}(t_{M+1}-t_M,\cdot)$ be a random function
distributed like ${y}(t_{M+1}-t_M,\cdot)$ but independent 
of $\{y(t,r)\}$.

Let $f$ be a bounded Lipschitz
 continuous function on $\Rb$, with Lipschitz constant $C_f$. 
The goal is to show
that the top line \eqref{yn-temp-7-0} below vanishes
 as $n\to\infty$. Add and
subtract terms to decompose \eqref{yn-temp-7-0} into three
differences: 
\begin{align}
&\Ef f\Bigl(\;\sum_{\substack{1\leq i\leq M+1\\[2pt] 1\leq j\leq N}}
 \theta_{i,j}
y_n(t_i,r_j) \Bigr) - Ef\Bigl(\;
\sum_{\substack{1\leq i\leq M+1\\[2pt] 1\leq j\leq N}} \theta_{i,j}
y(t_i,r_j) \Bigr) \label{yn-temp-7-0}\\
&=\; \biggl\{\; \Ef f\Bigl(\;\sum_{\substack{1\leq i\leq M+1\\[2pt] 1\leq j\leq N}}
 \theta_{i,j}
y_n(t_i,r_j) \Bigr)\nn\\
&\qquad\qquad  -
\Ef f\Bigl(\; \sum_{\substack{1\leq i\leq M\\[2pt] 1\leq k\leq K}}
 \rho^\om_{n,i,k}\,
y_n(t_i,v_k) + \sum_{1\leq j\leq N}
 \theta_{M+1,j}\,\tilde{y}_n(s_n,r_j)\Bigr)
\label{yn-temp-7-1}\;\biggr\}\\
&+\;\biggl\{\;  \Ef f\Bigl(\; \sum_{\substack{1\leq i\leq M\\[2pt] 1\leq k\leq K}}
 \rho^\om_{n,i,k}\,
y_n(t_i,v_k) + \sum_{1\leq j\leq N}
 \theta_{M+1,j}\,\tilde{y}_n(s_n,r_j)\Bigr)
\nn\\ 
&\qquad\qquad -
E f\Bigl(\; \sum_{\substack{1\leq i\leq M\\[2pt] 1\leq k\leq K}} \rho_{i,k}\,
y(t_i,v_k) + \sum_{1\leq j\leq N} \theta_{M+1,j}\,
\tilde{y}(t_{M+1}-t_M,r_j)\Bigr)
\label{yn-temp-7-2}\;\biggr\}\\
&+\;\biggl\{\; 
 E f\Bigl( \;\sum_{\substack{1\leq i\leq M\\[2pt] 1\leq k\leq K}} \rho_{i,k}\,
y(t_i,v_k) + \sum_{1\leq j\leq N}
 \theta_{M+1,j}\,\tilde{y}(t_{M+1}-t_M,r_j)\Bigr)\nn\\
&\qquad\qquad
- Ef\Bigl(\;\sum_{\substack{1\leq i\leq M+1\\[2pt] 1\leq j\leq N}} \theta_{i,j}
y(t_i,r_j) \Bigr)\;\biggr\}.
\label{yn-temp-7-3}
\end{align}
The remainder of the proof consists in treating the three
differences of expectations \eqref{yn-temp-7-1}--\eqref{yn-temp-7-3}. 

By the Lipschitz assumption and \eqref{yn-temp-6}, the difference
\eqref{yn-temp-7-1} is bounded  by 
\[
C_f\Ef\lvert R_n(a)\rvert +\Oc(n^{-1/4}).
\]
We need to bound $R_n(a)$. Recall that 
 $\gamma$ is an $\Nc(0, \sigma_a^2(t_{M+1}-t_M))$-distribution.

\begin{lemma} There exist  constants $C_1$ and $a_0$  such that,  
if  $a>a_0$, then for any partition $\{u_\ell\}$ of
$[-a,a]$  with mesh $\Delta$, 
and for any $1\leq j\leq N$, 
\[
\limsup_{n\to\infty} \Ef\lvert R_{n,j}(a)\rvert \leq C_1\bigl(
\sqrt{\Delta}+ \gamma(-\infty,-a/2)+\gamma(a/2,\infty) \bigr).
\]
\label{R(a)-lm}
\end{lemma}

We postpone the proof of Lemma \ref{R(a)-lm}. 
From this lemma, given $\e>0$, we can
choose first $a$ large enough and then $\Delta$ small enough 
so that 
\[
\limsup_{n\to\infty}\, [\text{ difference \eqref{yn-temp-7-1} }]
\leq \e/2.
\]
  
Difference \eqref{yn-temp-7-2} vanishes as $n\to\infty$, 
due to the induction assumption \eqref{yn-ind-M},   the limits 
 $\rho^\om_{n,i,k}\to\rho_{i,k}$ 
in probability, and the next lemma.  Notice that we are
not trying the invoke the induction assumption 
\eqref{yn-ind-M} for $M+1$ time points 
$\{t_1,\dotsc, t_M, s_n\}$. Instead, the induction
assumption is applied  to the first sum inside
$f$ in \eqref{yn-temp-7-2}. To the second sum apply
Lemma \ref{yn-lemma-1}, noting that 
$s_n\to t_{M+1}-t_M$. The two sums 
are  independent of each
other, as already observed after \eqref{yn-temp-6}, so 
they converge jointly.  This point is made precise in the next
lemma. 

\begin{lemma} Fix a positive integer $k$. 
For each $n$, let $V_n=(V_n^1,\dotsc,V^k_n)$, 
$X_n=(X_n^1,\dotsc,X^k_n)$, and  $\zeta_n$ be random 
variables on a common probability space.  Assume that 
 $X_n$ and $\zeta_n$ are independent of each other for each $n$.  
Let $v$ be a 
constant $k$-vector, $X$ another random $k$-vector, 
and $\zeta$ a random variable. Assume the weak limits
 $V_n\to v$, $X_n\to X$, and $\zeta_n\to \zeta$ hold marginally.
Then we have the weak limit 
\[
V_n\cdot X_n +\zeta_n\; \to\; v\cdot X +\zeta
\]
where the  $X$ and $\zeta$ on the right are independent. 
\label{weak-c-1-lm}
\end{lemma}

To prove this lemma, write
\[
V_n\cdot X_n +\zeta_n = (V_n-v)\cdot X_n + v\cdot X_n +\zeta_n
\]
and note that since $V_n\to v$ in probability, 
 tightness of $\{X_n\}$ implies that 
$(V_n-v)\cdot X_n \to 0$ in probability. As mentioned,
it applies to show that 
\[
\lim_{n\to\infty}\, [\text{ difference \eqref{yn-temp-7-2} }]=0.
\]

It remains to examine the  difference \eqref{yn-temp-7-3}.
From a consideration of how the coefficients $\rho^\om_{n,i,k}$
in \eqref{yn-temp-6} arise  and from the limit \eqref{q-clt-2},
\begin{align*}
&\sum_{\substack{1\leq i\leq M\\[2pt] 1\leq k\leq K}} \rho_{i,k}\,
y(t_i,v_k) + \sum_{1\leq j\leq N} \theta_{M+1,j}\,
\tilde{y}(t_{M+1}-t_M,r_j)
= 
\sum_{\substack{1\leq i\leq M\\[2pt] 1\leq j\leq N}} \theta_{i,j}
y(t_i,r_j)\\
&+ \sum_{1\leq j\leq N} \theta_{M+1,j}
\Bigl(\; \sum_{\ell=0}^{L-1} \gamma(u_\ell-r_j, u_{\ell+1}-r_j]
  y(t_M, u_\ell) + \tilde{y}(t_{M+1}-t_M,r_j) \Bigr)
\end{align*}
The first sum after the equality sign matches  all but the  
$(i=M+1)$-terms in the last sum in \eqref{yn-temp-7-3}. 
By virtue of the Markov property in 
\eqref{y-markov-1} we can represent the variables 
$y(t_{M+1},r_j)$ in the last sum in \eqref{yn-temp-7-3} by 
\[
y(t_{M+1},r_j)= \int_\Rb \varphi_{\sigma_a^2(t_{M+1}-t_M)}(u-r_j)
y(t_{M},u) \,d u
+ \tilde{y}(t_{M+1}-t_M, r_j).
\]
Then by the Lipschitz property of $f$ it  suffices to show 
that, for each $1\leq j\leq N$,  the expectation
\[
E\Bigl\lvert  
\int_\Rb \varphi_{\sigma_a^2(t_{M+1}-t_M)}(u-r_j)y(t_{M},u) \,d u
-
\sum_{\ell=0}^{L-1} \gamma(u_\ell-r_j, u_{\ell+1}-r_j]
  y(t_M, u_\ell) \Bigr\rvert
\]
can be made small by choice of $a>0$ and the partition $\{u_\ell\}$.
This follows from the moment bounds \eqref{y-mom-1}
 on the increments of 
the $y$-process and we omit the details. We have shown
that if  $a$ is large enough and then $\Delta$ small enough, 
\[
\limsup_{n\to\infty}\, [\text{ difference \eqref{yn-temp-7-3} }]
\leq \e/2.
\] 

To summarize, given
bounded Lipschitz $f$ and  $\e>0$, by choosing $a>0$ large
enough and the partition $\{u_\ell\}$ of $[-a,a]$ fine
enough,  
\[
\limsup_{n\to\infty} \biggl\lvert\;
\Ef f\Bigl(\;\sum_{\substack{1\leq i\leq M+1\\[2pt] 1\leq j\leq N}}
 \theta_{i,j}
y_n(t_i,r_j) \Bigr) - Ef\Bigl(\;
\sum_{\substack{1\leq i\leq M+1\\[2pt] 1\leq j\leq N}} \theta_{i,j}
y(t_i,r_j) \Bigr)\;\biggr\rvert \leq \e.
\]
This completes the proof of the induction step 
and thereby the proof of Theorem \ref{yn-thm-1}. 
\end{proof} 

It remains  to verify the lemmas that were used along the way. 

\begin{proof}[Proof of Lemma \ref{R(a)-lm}]
We begin with a calculation. Here it is convenient to use
the space-time walk $\Xbar^{x,m}_k=(X^{x,m}_k, m-k)$. 
First observe that 
\be
\begin{aligned}
&E^\om(X^{x,m}_n)-x-nV=
 \sum_{k=0}^{n-1} E^{\om}\bigl[
X^{x,m}_{k+1}-X^{x,m}_k-V \bigr]\\
&= \sum_{k=0}^{n-1} E^{\om}\bigl[
E^{ T_{\{\Xbar^{x,m}_k\}}\om} (X^{0,0}_1-V)\bigr]
=\sum_{k=0}^{n-1} E^{\om} g( T_{\Xbar^{x,m}_k}\om).
\end{aligned}
\label{rw-calc-1} 
\ee
From this, for $x,y\in\Zb$, 
\begin{align*}
&\Ef \Bigl[
\bigl(\{E^\om(X^{x,n}_n)-x\}-  \{E^\om(X^{y,n}_n)-y\}\bigr)^2
\Bigr]\\
&= \Ef \sum_{k=0}^{n-1} \Bigl( E^{\om} 
\bigl[g(T_{\Xbar^{x,n}_k}\om) -g(T_{\Xbar^{y,n}_k}\om) \bigr] \Bigr)^2\\
&\qquad +2  \sum_{0\leq k <\ell<n} \Ef  E^{\om} 
\bigl[g(T_{\Xbar^{x,n}_k}\om) -g(T_{\Xbar^{y,n}_k}\om) \bigr] 
E^{\om} \bigl[g(T_{\Xbar^{x,n}_\ell}\om) -g(T_{\Xbar^{y,n}_\ell}\om) \bigr]\\ 
&\qquad\qquad\qquad\qquad \text{(the cross terms for $k<\ell$ vanish)}\\
&= \Ef \sum_{k=0}^{n-1} \Bigl( \sum_{z,w\in\Zb^2} P^{\om}\{\Xbar^{x,n}_k=z\}
P^{\om}\{\Xbar^{y,n}_k=w\}   
\bigl[g(T_{z}\om) -g(T_{w}\om) \bigr] \Bigr)^2\\
&= \Ef \sum_{k=0}^{n-1}  \sum_{z,w,u,v\in\Zb^2} 
P^{\om}\{\Xbar^{x,n}_k=z\}P^{\om}\{\Xbar^{y,n}_k=w\}   
P^{\om}\{\Xbar^{x,n}_k=u\}P^{\om}\{\Xbar^{y,n}_k=v\}   \\
&\qquad \times \Bigl( g(T_{z}\om)g(T_{u}\om) -g(T_{w}\om)g(T_{u}\om)
-g(T_{z}\om)g(T_{v}\om)+ g(T_{w}\om)g(T_{v}\om) \Bigr)\\
&\qquad\qquad\qquad\qquad
\text{(by independence $\Ef g(T_{z}\om)g(T_{u}\om)
= \sigma_D^2 \ind_{\{z=u\}}$)}\\
&=\sigma_D^2  \sum_{k=0}^{n-1} \Bigl( P\{X^{x,n}_k=\Xtil^{x,n}_k\}
- 2P\{X^{x,n}_k=\Xtil^{y,n}_k\} +P\{X^{y,n}_k=\Xtil^{y,n}_k\}  \Bigr)\\
&=2\sigma_D^2  \sum_{k=0}^{n-1} \Bigl( P_0\{Y_k=0\}
- P_{x-y}\{Y_k=0\} \Bigr) \\
&=2\sigma_D^2 \bigl( G_{n-1}(0,0)-G_{n-1}(x-y,0)\bigr).
\end{align*}
On the last three lines above, as elsewhere in the paper, we 
used these conventions: $X_k$ and $\Xtil_k$ denote walks that are
independent in a common environment $\om$, 
$Y_k=X_k-\Xtil_k$ is the difference walk, and $G_n(x,y)$ 
the Green function   of $Y_k$. By Lemma \ref{Gn-lm-2} 
 we get the inequality 
\be
\Ef \Bigl[
\bigl(\{E^\om(X^{x,n}_n)-x\}-  \{E^\om(X^{y,n}_n)-y\}\bigr)^2
\Bigr]
\leq C\abs{x-y}
\label{ineq-R-temp-1}
\ee
valid for all $n$ and all $x,y\in\Zb$. 

Turning to $R_{n,j}(a)$ defined in 
\eqref{def-R(a)-1}--\eqref{def-R(a)-2}, and 
utilizing independence, 
\begin{align*}
&\Ef\lvert R_{n,j}(a)\rvert  \leq 
n^{-1/4}\sum_{z=\tfl{-a\sqrt{n}}+1}^{\tfl{a\sqrt{n}}}
\Ef[p^\om_{n,j}(z)] 
\Bigl(  \Ef \Bigl[ \bigl( \{
E^\om(X^{\tfl{nt_Mb} +  z, \tfl{nt_M}}_{\tfl{nt_M}})
-z\}  \\
& \qquad\qquad\qquad\qquad -\;  \{
E^\om(X^{\tfl{nt_Mb}+  \tfl{u(z)\sqrt{n}},
 \tfl{nt_M}}_{\tfl{nt_M}})
-\tfl{u(z)\sqrt{n}}\} \,\bigr)^2 \Bigr]\;\Bigr)^{1/2} 
\\
&+ n^{-1/4}\sum_{\substack{z\leq -a\sqrt{n}\\ z> a\sqrt{n}}}
\Ef[p^\om_{n,j}(z)] 
\Bigl( \Ef \bigl[  \bigl(
E^\om(X^{\tfl{nt_Mb} +  z, \tfl{nt_M}}_{\tfl{nt_M}})
-E(X^{\tfl{nt_Mb} +  z, \tfl{nt_M}}_{\tfl{nt_M}})
    \bigr)^2\,\bigr]\Bigr)^{1/2}\\
&\qquad\qquad + n^{-1/4}\sum_{\substack{z\leq -a\sqrt{n}\\z> a\sqrt{n}}}
\Ef[p^\om_{n,j}(z)]\cdot  
 \bigl\lvert \,
E(X^{\tfl{nt_Mb} +  z, \tfl{nt_M}}_{\tfl{nt_M}})
-z\,\bigr\rvert \\
&\leq  Cn^{-1/4} 
 \max_{-a\sqrt{n} <z\leq   a\sqrt{n}} \,\lvert z-
\tfl{u(z)\sqrt{n}} \,\rvert^{1/2} \\ 
&\quad +CP\bigl\{ \,\bigl\lvert X^{\tfl{nt_{M+1}b}+\tfl{r_j\sqrt{n}}, 
\tfl{nt_{M+1}}}_{\tfl{nt_{M+1}}-\tfl{nt_{M}}}-\tfl{nt_Mb}\bigr\rvert
\geq a\sqrt{n} \,\bigr\}
  +C n^{-1/4}. 
\end{align*}
For the last inequality above we used \eqref{ineq-R-temp-1},
bound \eqref{qm-var-1} on the variance of the quenched mean,
and then 
\[
E(X^{\tfl{nt_Mb} +  z, \tfl{nt_M}}_{\tfl{nt_M}})
-z = \tfl{nt_Mb}+ \tfl{nt_M}V = \tfl{nt_Mb}- \tfl{nt_M}b=\Oc(1).
\]
By the choice of $u(z)$, and
 by the central limit theorem  if $a> 2\lvert r_j\rvert$, 
the limit of the bound
on $\Ef\lvert R_{n,j}(a)\rvert$ 
as $n\to\infty$ is  $C(\sqrt\Delta +\gamma(-\infty, -a/2)
+\gamma(a/2,\infty))$. This completes the proof of 
Lemma \ref{R(a)-lm}. 
\end{proof}

\begin{proof}[Proof of Lemma \ref{yn-lemma-1}]
 We drop the subscript 
from $t_n$ and write simply $t$.  For the main part of the 
proof the only relevant property is that $nt_n=O(n)$. 
We point this out after the preliminaries. 

We show convergence of the linear combination
$\sum_{i=1}^N \theta_iy_n(t,r_i)$ for an arbitrary but fixed 
$N$-vector $\theta=(\theta_1,\dotsc,\theta_N)$. This in turn
will come from a martingale central limit theorem. 
For this proof abbreviate 
$X^i_{k}= X^{\tfl{ntb}+\tfl{r_i\sqrt{n}}, \tfl{nt}}_{k}$. 
For $1\leq k\leq \tfl{nt}$ define
\[
z_{n,k}=n^{-1/4}\sum_{i=1}^N \theta_i E^\om g(T_{\Xbar^i_{k-1}}\om)
\]
so that by \eqref{rw-calc-1} 
\[
\sum_{k=1}^{\tfl{nt}} z_{n,k}= \sum_{i=1}^N \theta_iy_n(t,r_i)
+\Oc(n^{-1/4}).
\]
The error is deterministic and comes from the 
discrepancy \eqref{yn-center}
in the centering. It vanishes in the limit and so can be ignored. 

A probability of the type 
$P^\om(X^i_{k-1}=y)$ is a function of the 
environments  \[\{\ombar_{j}:  \tfl{nt}-k+2\leq j\leq \tfl{nt}\}\]
while $g(T_{y,s}\om)$ is a function of 
$\ombar_{s}$.
For a fixed $n$, $\{z_{n,k}: 1\leq k\leq \tfl{nt}\}$ are martingale 
differences  with respect to the filtration 
\[
\Uc_{n,k}=\sigma\{ \ombar_{j}:  \tfl{nt}-k+1\leq j\leq \tfl{nt} \} 
\qquad (1\leq k\leq \tfl{nt})
\]
with $\Uc_{n,0}$ equal to the  trivial $\sigma$-algebra. 
The goal is to show that $\sum_{k=1}^{\tfl{nt}} z_{n,k}$ 
converges to a centered Gaussian with variance 
$\sum_{1\leq i,j \leq N} \theta_i\theta_j \Gamma_q((t,r_i),(t,r_j))$. 
By the Lindeberg-Feller Theorem for martingales, 
 it suffices to check  that 
\be
\sum_{k=1}^{\tfl{nt}} \Ef[z_{n,k}^2\,\vert\, \Uc_{n,k-1}]
\longrightarrow 
\sum_{1\leq i,j \leq N} \theta_i\theta_j \Gamma_q((t,r_i),(t,r_j))
\label{check-1}
\ee
and 
\be
\sum_{k=1}^{\tfl{nt}} \Ef[z_{n,k}^2\ind\{\lvert z_{n,k}\rvert\geq \e\}
\,\vert\, \Uc_{n,k-1}]
\longrightarrow 0
\label{check-2}
\ee
in probability, as $n\to \infty$,  for every $\e>0$. 
Condition \eqref{check-2} is trivially satisfied because
$\lvert z_{n,k}\rvert \leq Cn^{-1/4}$ by the 
boundedness of $g$. 

The main part of the proof consists 
of checking \eqref{check-1}. 
This argument is a generalization of the proof of
 \cite[Theorem 4.1]{ferr-font-rap} where it was done for a nearest-neighbor
walk.  We follow their reasoning for the first part of the proof. 
 Since 
$\sigma_D^2=\Ef[g^2]$ and since conditioning $z_{n,k}^2$ on
$\Uc_{n,k-1}$ entails integrating out the environments
$\ombar_{\tfl{nt}-k+1}$, one can derive 
\[
\sum_{k=1}^{\tfl{nt}} \Ef[z_{n,k}^2\,\vert\, \Uc_{n,k-1}]
=\sigma_D^2\sum_{1\leq i,j \leq N} \theta_i\theta_j\,
n^{-1/2} \sum_{k=0}^{\tfl{nt}-1}
P^\om(X^i_k=\Xtil^j_k)
\]
where $X^i_k$ and $\Xtil^j_k$ are two walks independent 
under the common environment $\om$, started at 
$(\tfl{ntb}+\tfl{r_i\sqrt{n}}, \tfl{nt})$ and 
$(\tfl{ntb}+\tfl{r_j\sqrt{n}}, \tfl{nt})$.  

By \eqref{G-lim1} 
\be
\sigma_D^2 n^{-1/2} \sum_{k=0}^{\tfl{nt}-1}
P(X^i_k=\Xtil^j_k)\longrightarrow \Gamma_q((t,r_i),(t,r_j)).
\label{G-q-appear}
\ee
This limit holds if instead of a fixed $t$ on the left 
we have a sequence $t_n\to t$. 
Consequently we will have proved \eqref{check-1} if we show,
for each fixed pair $(i,j)$,  that 
\be
n^{-1/2} \sum_{k=0}^{\tfl{nt}-1} \bigl(P^\om\{X^i_k=\Xtil^j_k\}-
P\{X^i_k=\Xtil^j_k\}\bigr)\longrightarrow 0
\label{check-1-1}
\ee
in $\Pf$-probability. 
For the above statement the behavior of $t$ is immaterial as
long as it stays bounded as $n\to\infty$.

Rewrite the expression in \eqref{check-1-1}  as 
\begin{align*}
&n^{-1/2} \sum_{k=0}^{\tfl{nt}-1} \bigl(
P\{X^i_k=\Xtil^j_k\,\vert\,\Uc_{n,k}\}-
P\{X^i_k=\Xtil^j_k\,\vert\,\Uc_{n,0}\}\bigr)\\
&=n^{-1/2} \sum_{k=0}^{\tfl{nt}-1} \sum_{\ell=0}^{k-1} \bigl(
P\{X^i_k=\Xtil^j_k\,\vert\,\Uc_{n,\ell+1}\}-
P\{X^i_k=\Xtil^j_k\,\vert\,\Uc_{n,\ell}\}\bigr)\\
&=n^{-1/2} \sum_{\ell=0}^{\tfl{nt}-1}
 \sum_{k=\ell+1}^{\tfl{nt}-1}  \bigl(
P\{X^i_k=\Xtil^j_k\,\vert\,\Uc_{n,\ell+1}\}-
P\{X^i_k=\Xtil^j_k\,\vert\,\Uc_{n,\ell}\}\bigr)\\
&\equiv n^{-1/2} \sum_{\ell=0}^{\tfl{nt}-1} R_\ell
\end{align*}
where the last line defines $R_\ell$. 
Check that $\Ef R_\ell R_m =0$ for $\ell\neq m$. 
Thus it is convenient to verify 
our goal \eqref{check-1-1} by checking $L^2$ convergence,
in other words by showing 
\begin{align}
 &n^{-1} \sum_{\ell=0}^{\tfl{nt}-1} \Ef[R_\ell^2]\nn\\
&=n^{-1} \sum_{\ell=0}^{\tfl{nt}-1}
\Ef\biggl[\;\biggl\{ \;  \sum_{k=\ell+1}^{\tfl{nt}-1}  \bigl(
P\{X^i_k=\Xtil^j_k\,\vert\,\Uc_{n,\ell+1}\}-
P\{X^i_k=\Xtil^j_k\,\vert\,\Uc_{n,\ell}\}\bigr)\biggr\}^2\,\biggr]
\label{check-1-2}\\
&\longrightarrow 0. \nn
\end{align}

For the moment we work on a single
 term inside the braces in \eqref{check-1-2}, for 
a fixed pair $k>\ell$.
Write $Y_m=X^i_m-\Xtil^j_m$ for the difference
walk. By the Markov property of the walks [recall \eqref{X-tr-pr}]
 we can write 
\begin{align*}
&P\{X^i_k=\Xtil^j_k\,\vert\,\Uc_{n,\ell+1}\}
=
\sum_{x,\xtil,y,\ytil\in\Zb} P^\om\{ X^i_\ell=x,\Xtil^j_\ell=\xtil\} \\
&\qquad\qquad \times
 u^\om_{\tfl{nt}-\ell}(x,y-x)u^\om_{\tfl{nt}-\ell}(\xtil,\ytil-\xtil)
P(Y_k=0\,\vert\, Y_{\ell+1}=y-\ytil)
\end{align*}
and similarly for the other conditional probability
\begin{align*}
&P\{X^i_k=\Xtil^j_k\,\vert\,\Uc_{n,\ell}\}
=
\sum_{x,\xtil,y,\ytil\in\Zb} P^\om\{ X^i_\ell=x,\Xtil^j_\ell=\xtil\} \\
&\qquad\qquad \times
 \Ef[u^\om_{\tfl{nt}-\ell}(x,y-x)u^\om_{\tfl{nt}-\ell}(\xtil,\ytil-\xtil)]
P(Y_k=0\,\vert\, Y_{\ell+1}=y-\ytil).
\end{align*}
Introduce the transition probability $q(x,y)$ of the $Y$-walk.
Combine the above decompositions
 to express the $(k,\ell)$ term inside the braces in
\eqref{check-1-2}  as
\begin{align}
&P\{X^i_k=\Xtil^j_k\,\vert\,\Uc_{n,\ell+1}\}-
P\{X^i_k=\Xtil^j_k\,\vert\,\Uc_{n,\ell}\} \nn\\
&=\sum_{x,\xtil,y,\ytil\in\Zb} P^\om\{ X^i_\ell=x,\Xtil^j_\ell=\xtil\} 
q^{k-\ell-1}(y-\ytil,0)  \nn \\
&\qquad\qquad \times
\bigl( u^\om_{\tfl{nt}-\ell}(x,y-x)u^\om_{\tfl{nt}-\ell}(\xtil,\ytil-\xtil)
- \Ef[u^\om_{\tfl{nt}-\ell}(x,y-x)u^\om_{\tfl{nt}-\ell}(\xtil,\ytil-\xtil)]
\,\bigr)\nn\\
&=\sum_{x,\xtil} P^\om\{ X^i_\ell=x,\Xtil^j_\ell=\xtil\} 
\sum_{\substack{ z,w:\; -\rg\leq w\leq \rg\\[2pt] -\rg\leq w-z\leq \rg}}
  q^{k-\ell-1}(x-\xtil+z,0) \nn\\
&\qquad\qquad \times 
\Bigl( u^\om_{\tfl{nt}-\ell}(x,w)
u^\om_{\tfl{nt}-\ell}(\xtil,w-z)\nn\\
&\qquad\qquad \qquad\qquad\qquad - \Ef[u^\om_{\tfl{nt}-\ell}(x,w)
 u^\om_{\tfl{nt}-\ell}(\xtil,w-z)]
\,\Bigr). \label{check-1-3a}
\end{align}
The last sum above uses  the finite
range $\rg$ of the jump probabilities. 
Introduce the quantities
\[
\rho^\om_\ell(x,x+m)=\sum_{y: y\leq m} u^\om_{\tfl{nt}-\ell}(x,y)
=\sum_{y=-\rg}^{m} u^\om_{\tfl{nt}-\ell}(x,y)
\]
and 
\begin{align*}
\zeta^\om_\ell(x,\xtil,z,w) &=
  \rho^\om_{\ell}(x,x+w)
u^\om_{\tfl{nt}-\ell}(\xtil,w-z)\\
& \qquad\qquad\qquad - \Ef\bigl[
 \rho^\om_{\ell}(x,x+w) 
 u^\om_{\tfl{nt}-\ell}(\xtil,w-z)\bigr]. 
\end{align*}
Fix $(x,\xtil)$, consider the sum over $z$ and $w$ on line
 \eqref{check-1-3a}, and 
continue with  a ``summation by parts'' step:   
\begin{align*}
&\sum_{\substack{ z,w:\; -\rg\leq w\leq \rg\\[2pt] -\rg\leq w-z\leq \rg}}
 q^{k-\ell-1}(x-\xtil+z,0) 
  \Bigl( u^\om_{\tfl{nt}-\ell}(x,w)
u^\om_{\tfl{nt}-\ell}(\xtil,w-z)\\
&\qquad\qquad\qquad \qquad\qquad\qquad - \Ef[u^\om_{\tfl{nt}-\ell}(x,w)
 u^\om_{\tfl{nt}-\ell}(\xtil,w-z)]
\,\Bigr)\\
&=\sum_{\substack{ z,w:\; -\rg\leq w\leq \rg\\[2pt] -\rg\leq w-z\leq \rg}}
 q^{k-\ell-1}(x-\xtil+z,0) 
 \bigl( \zeta^\om_{\ell}(x,\xtil,z, w) - 
\zeta^\om_{\ell}(x,\xtil,z-1, w-1)  \bigr) \\
&=\sum_{\substack{ z,w:\; -\rg\leq w\leq \rg\\[2pt] -\rg\leq w-z\leq \rg}}
 \Bigl( q^{k-\ell-1}(x-\xtil+z,0)
- q^{k-\ell-1}(x-\xtil+z+1,0)\Bigr) \zeta^\om_{\ell}(x,\xtil,z, w) \\
&\qquad\qquad\qquad  +\sum_{z=0}^{2\rg}
q^{k-\ell-1}(x-\xtil+z+1,0) \zeta^\om_{\ell}(x,\xtil,z, \rg) \\
&\qquad\qquad\qquad - \sum_{z=-2\rg-1}^{-1}
q^{k-\ell-1}(x-\xtil+z+1,0) \zeta^\om_{\ell}(x,\xtil,z, -\rg-1). 
\end{align*}
By definition of the 
range $\rg$, the last sum above vanishes because 
 $\zeta^\om_{\ell}(x,\xtil,z, -\rg-1)=0$. 
Take this into consideration, 
substitute the last form above into \eqref{check-1-3a} 
and sum over  $k=\ell+1,\dotsc, \tfl{nt}-1$. 
Define the quantity
\be
A_{\ell,n}(x)=\sum_{k=\ell+1}^{\tfl{nt}-1} 
\bigl( q^{k-\ell-1}(x,0)
- q^{k-\ell-1}(x+1,0) \bigr).  
\label{def-Aellnx}
\ee
Then 
the expression in braces in \eqref{check-1-2} is represented as 
\begin{align}
&R_\ell = \sum_{k=\ell+1}^{\tfl{nt}-1}  \bigl(
P\{X^i_k=\Xtil^j_k\,\vert\,\Uc_{n,\ell+1}\}-
P\{X^i_k=\Xtil^j_k\,\vert\,\Uc_{n,\ell}\}\bigr)\nn\\
&=\sum_{x,\xtil}  P^\om\{ X^i_\ell=x,\Xtil^j_\ell=\xtil\}
\sum_{\substack{ z,w:\; -\rg\leq w\leq \rg\\[2pt] -\rg\leq w-z\leq \rg}}
A_{\ell,n}(x-\xtil+z) \zeta^\om_\ell(x,\xtil,z,w)
\label{check-1-4}\\
 +&\sum_{x,\xtil}  P^\om\{ X^i_\ell=x,\Xtil^j_\ell=\xtil\}
\sum_{z=0}^{2\rg}
\sum_{k=\ell+1}^{\tfl{nt}-1}  
q^{k-\ell-1}(x-\xtil+z+1,0) \zeta^\om_{\ell}(x,\xtil,z, \rg)
\label{check-1-4a}\\
&\equiv R_{\ell,1}+R_{\ell,2}\nn
\end{align}
where  $R_{\ell,1}$ and $R_{\ell,2}$ denote the sums on lines
\eqref{check-1-4} and \eqref{check-1-4a}. 

Recall from \eqref{check-1-2} that our goal was to show
that $n^{-1}\sum_{\ell=0}^{\tfl{nt}-1} \Ef R_\ell^2\to 0$
 as $n\to\infty$.  We show this separately for
 $R_{\ell,1}$ and $R_{\ell,2}$. 

As a function of $\om$, $\zeta^\om_{\ell}(\dotsm)$ is a 
function of $\ombar_{\tfl{nt}-\ell}$ and hence 
independent of the probabilities on line 
\eqref{check-1-4}. 
Thus we get 
\begin{align}
\Ef[R_{\ell,1}^2] &=
\sum_{x,\xtil,x',\xtil'}
\Ef\bigl[ P^\om\{ X^i_\ell=x,\Xtil^j_\ell=\xtil\}
P^\om\{ X^i_\ell=x',\Xtil^j_\ell=\xtil'\}\bigr] 
 \nn\\
&\qquad \times 
\sum_{\substack{ -\rg\leq w\leq \rg\\[2pt] -\rg\leq w-z\leq \rg}}
\sum_{\substack{  -\rg\leq w'\leq \rg\\[2pt] -\rg\leq w'-z'\leq \rg}}
A_{\ell,n}(x-\xtil+z)A_{\ell,n}(x'-\xtil'+z') \nn\\
&\qquad\qquad\qquad\qquad \times 
\Ef\bigl[ \zeta^\om_\ell(x,\xtil,z,w)
\zeta^\om_\ell(x',\xtil',z',w')  \bigr]
\label{check-1-5}
\end{align}
Lemma \ref{Gn-lm-2} implies that $A_{\ell, n}(x)$ 
is uniformly bounded over $(\ell,n,x)$. 
Random variable $\zeta^\om_\ell(x,\xtil,z,w)$ is mean zero
and a function of the environments 
$\{ \om_{x,\tfl{nt}-\ell}, \om_{\xtil,\tfl{nt}-\ell} \}$.
Consequently the last expectation  on line \eqref{check-1-5}
vanishes unless $\{x,\xtil\}\cap\{x',\xtil'\}\neq\emptyset$. 
The sums over $z,w,z',w'$ contribute
a constant because of their bounded range. Taking all these
into consideration, we obtain the bound 
\be
\Ef[R_{\ell,1}^2] \leq C\Bigl(
P\{X^i_\ell=\Xtil^i_\ell\}
+P\{X^i_\ell=\Xtil^j_\ell\}+P\{X^j_\ell=\Xtil^j_\ell\}\Bigr).
\label{check-1-6}
\ee
By \eqref{G-lim1} we get the bound
\be
n^{-1}\sum_{\ell=0}^{\tfl{nt}-1} \Ef [R_{\ell,1}^2]
\leq Cn^{-1/2} 
\label{check-1-7}
\ee
which vanishes as $n\to\infty$. 

For the remaining sum $R_{\ell,2}$ 
observe first that 
\be
\zeta^\om_{\ell}(x,\xtil,z, \rg)= 
u^\om_{\tfl{nt}-\ell}(\xtil,\rg-z) 
- \Ef  u^\om_{\tfl{nt}-\ell}(\xtil,\rg-z).
\ee
 Summed over $0\leq z\leq 2\rg$ this vanishes, so we can
start by rewriting as follows: 
\begin{align*}
&R_{\ell,2}= 
\sum_{x,\xtil}  P^\om\{ X^i_\ell=x,\Xtil^j_\ell=\xtil\}\nn\\
&\qquad    \times 
\sum_{z=0}^{2\rg} \sum_{k=\ell+1}^{\tfl{nt}-1}  
\bigl( q^{k-\ell-1}(x-\xtil+z+1,0)-q^{k-\ell-1}(x-\xtil,0)\bigr) 
\zeta^\om_{\ell}(x,\xtil,z, \rg)\\
&=- \sum_{x,\xtil}  P^\om\{ X^i_\ell=x,\Xtil^j_\ell=\xtil\}
\sum_{z=0}^{2\rg} \sum_{m=0}^z A_{\ell,n}(x-\xtil+m,0) 
\zeta^\om_{\ell}(x,\xtil,z, \rg)\\
&=- \sum_{x,\xtil}  P^\om\{ X^i_\ell=x,\Xtil^j_\ell=\xtil\}
\sum_{m=0}^{2\rg}  A_{\ell,n}(x-\xtil+m,0) 
\rhobar^\om_{\ell}(\xtil,\xtil+ \rg-m)
\end{align*}
where we abbreviated on the last line
\[
\rhobar^\om_{\ell}(\xtil,\xtil+ \rg-m)
=\rho^\om_{\ell}(\xtil,\xtil+ \rg-m)
- \Ef \rho^\om_{\ell}(\xtil,\xtil+ \rg-m).
\]
Square the last representation for $R_{\ell,2}$, 
take $\Ef$-expectation, and note that
\[
\Ef\bigl[ \rhobar^\om_{\ell}(\xtil,\xtil+ \rg-m)
\rhobar^\om_{\ell}(\xtil',\xtil'+ \rg-m')\bigr]=0
\]
unless $\xtil=\xtil'$. Thus the reasoning applied to $R_{\ell,1}$
can be repeated, and we conclude that also 
 $n^{-1}\sum_{\ell=0}^{\tfl{nt}-1} \Ef R_{\ell,2}^2\to 0$.

To summarize,  we have verified 
\eqref{check-1-2}, thereby \eqref{check-1-1} and  
condition \eqref{check-1} for the martingale CLT. 
This completes the proof of Lemma \ref{yn-lemma-1}.
\end{proof}

\section{Proofs for forward walks in a random environment}
\label{fw-walk-pf-sec} 

\subsection{Proof of Theorem \ref{an-thm-1}}
The proof of Theorem \ref{an-thm-1} is organized in the same 
way as the proof of  Theorem \ref{yn-thm-1} so we restrict
ourselves to a few remarks. 
The Markov property reads now $(0\leq s<t, r\in\Rb)$:
\begin{align*}
&a_n(t,r)=a_n(s,r)+ \sum_{y\in\Zb} P^\om\bigl\{ 
Z^{\tfl{r\sqrt{n}},0}_{\tfl{ns}}=\tfl{r\sqrt{n}}+\tfl{nsV}+y\bigr\}\\
&\quad \times n^{-1/4}  \bigl\{ 
E^\om( Z^{\tfl{r\sqrt{n}}+\tfl{nsV}+y,\tfl{ns}}_{\tfl{nt}-\tfl{ns}} )
-\tfl{r\sqrt{n}}-y -\tfl{nt}V\bigr\}+n^{-1/4}\bigl(\tfl{ns}V-\tfl{nsV} \bigr). 
\end{align*}
This serves as the basis for the inductive proof along time 
levels, exactly as done in the argument following 
\eqref{yn-ind-M}. 

Lemma \ref{yn-lemma-1} about the convergence at a fixed $t$-level
applies to $a_n(t,\cdot)$ exactly as worded. This follows from 
noting that, up to a trivial difference from integer parts, 
the processes $a_n(t,\cdot)$ and $y_n(t,\cdot)$ are the same. 
Precisely, if $S$ denotes the $\Pf$-preserving
 transformation on $\Omega$
defined by $(S\om)_{x,\tau}=\om_{-\tfl{ntb}+x,\tfl{nt}-\tau}$,
then 
\[
E^{S\om}( X^{\tfl{ntb}+\tfl{r\sqrt{n}},\tfl{nt}}_{\tfl{nt}} )
-\tfl{r\sqrt{n}} =
E^\om( Z^{\tfl{r\sqrt{n}},0}_{\tfl{nt}} )-\tfl{r\sqrt{n}}
+\tfl{ntb}.
\]
The errors in the inductive argument are treated with the same 
arguments as used in Lemma \ref{R(a)-lm} to treat $R_{n,j}(a)$. 


\subsection{Proof of Corollary \ref{a-cor-1}}
We start with a moment bound that will give tightness of the 
processes.

\begin{lemma} There exists a constant $0<C<\infty$ such that,
for all $n\in\Nb$,  
\[
\Ef\bigl[ (E^\om(Z^{0,0}_n)-nV)^6\bigr]\leq Cn^{3/2}. 
\]
\label{qm-mom-6}
\end{lemma}

\begin{proof} From 
\[
\Ef\bigl[ \bigl( E^{\om} 
g(T_{\Zbar^{x,0}_n}\om) -E^{\om} g(T_{\Zbar^{0,0}_n}\om)
\bigr)^2 \bigr]
= 2\sigma_D^2\bigl( 
P[Y^0_n=0]- P[Y^x_n=0] \bigr) 
\]
we get 
\be
P[Y^x_n=0]\leq P[Y^0_n=0]\qquad\text{for all $n\geq 0$ and $x\in\Zb$.}
\label{Y-mon}
\ee
Abbreviate $\Zbar_n=\Zbar^{0,0}_n$ for this proof. 
 $E^\om(Z_n)-nV$ is  a mean-zero martingale
with increments $E^\om g(T_{\Zbar_k}\om)$
relative to the filtration 
$\Hc_n=\sigma\{\ombar_k: 0\leq k<n\}$. By the 
Burkholder-Davis-Gundy  inequality \cite{burk},
\begin{align*}
\Ef\bigl[ (E^\om(Z_n)-nV)^6\bigr]\leq  C
\Ef\biggl[ \Bigl(\;\sum_{k=0}^{n-1} \bigl[E^\om g(T_{\Zbar_k}\om)  
\bigr]^2\Bigr)^3\biggr].
\end{align*}
Expanding the cube yields four sums 
\begin{align*}
&C\sum_{0\leq k<n} \Ef \bigl[ \bigl(E^\om g(T_{\Zbar_k}\om)  
\bigr)^6\bigr]
+
C\sum_{0\leq k_1<k_2<n} \Ef \Bigl[ \bigl(E^\om g(T_{\Zbar_{k_1}}\om)  
\bigr)^2  \bigl(E^\om g(T_{\Zbar_{k_2}}\om)  
\bigr)^4\Bigr] \\
&\quad +
C\sum_{0\leq k_1<k_2<n} \Ef \Bigl[ \bigl(E^\om g(T_{\Zbar_{k_1}}\om)  
\bigr)^4  \bigl(E^\om g(T_{\Zbar_{k_2}}\om)  
\bigr)^2\Bigr] \\
&\quad\quad +
C\sum_{0\leq k_1<k_2<k_3<n} \Ef \Bigl[ \bigl(E^\om g(T_{\Zbar_{k_1}}\om)  
\bigr)^2  \bigl(E^\om g(T_{\Zbar_{k_2}}\om)  \bigr)^2
 \bigl(E^\om g(T_{\Zbar_{k_2}}\om)  \bigr)^2
\Bigr] 
\end{align*} 
with a constant $C$  that bounds  the number of 
arrangements of each type. 
Replacing some $g$-factors with constant upper bounds simplifies
the quantity  to this:
\begin{align*}
C\sum_{0\leq k<n} \Ef \bigl[ \bigl(E^\om g(T_{\Zbar_k}\om)  
\bigr)^2\bigr]
+
C\sum_{0\leq k_1<k_2<n} \Ef \Bigl[ \bigl(E^\om g(T_{\Zbar_{k_1}}\om)  
\bigr)^2  \bigl(E^\om g(T_{\Zbar_{k_2}}\om)  
\bigr)^2\Bigr] \\
\quad +
C\sum_{0\leq k_1<k_2<k_3<n} \Ef \Bigl[ \bigl(E^\om g(T_{\Zbar_{k_1}}\om)  
\bigr)^2  \bigl(E^\om g(T_{\Zbar_{k_2}}\om)  \bigr)^2
 \bigl(E^\om g(T_{\Zbar_{k_3}}\om)  \bigr)^2
\Bigr].  
\end{align*} 
The expression above is bounded by $C(n^{1/2}+n+n^{3/2})$. We show the 
 argument for the last sum of triple products. (Same reasoning
applies to the first two sums.)  It 
utilizes repeatedly independence, 
$\Ef g(T_{u}\om)g(T_{v}\om)
= \sigma_D^2 \ind_{\{u=v\}}$ for $u,v\in\Zb^2$, 
 and \eqref{Y-mon}. Fix $0\leq k_1<k_2<k_3<n$. 
Let $\Zbar'_k$ denote an  independent  copy 
of the walk $\Zbar_k$ in the same environment $\om$. 
\begin{align*}
&\Ef \Bigl[ \bigl(E^\om g(T_{\Zbar_{k_1}}\om)  
\bigr)^2  \bigl(E^\om g(T_{\Zbar_{k_2}}\om)  \bigr)^2
 \bigl(E^\om g(T_{\Zbar_{k_3}}\om)  \bigr)^2
\Bigr]\\
&=\Ef \biggl[ \bigl(E^\om g(T_{\Zbar_{k_1}}\om)  
\bigr)^2  \bigl(E^\om g(T_{\Zbar_{k_2}}\om)  \bigr)^2\\
&\qquad\qquad\qquad \times 
 \sum_{u,v\in\Zb^2} P^\om\{\Zbar_{k_3}=u\} P^\om\{\Zbar'_{k_3}=v\}
 \Ef\bigl\{ g(T_{u}\om) g(T_{v}\om)  \bigr\}\,
\biggr]\\
&=C\Ef \Bigl[ \bigl(E^\om g(T_{\Zbar_{k_1}}\om)  
\bigr)^2  \bigl(E^\om g(T_{\Zbar_{k_2}}\om)  \bigr)^2
 P^\om\{\Zbar_{k_3}=\Zbar'_{k_3}\}
 \Bigr]\\
&=C\Ef \biggl[ \bigl(E^\om g(T_{\Zbar_{k_1}}\om)  
\bigr)^2  \bigl(E^\om g(T_{\Zbar_{k_2}}\om)  \bigr)^2\\
&\qquad\qquad\qquad \times 
\sum_{u,v\in\Zb^2} P^\om\{\Zbar_{k_2+1}=u, \Zbar'_{k_2+1}=v\}
\biggr] \Ef P^\om\{ \Zbar^u_{k_3-k_2-1}=\Zbar^v_{k_3-k_2-1}\}  \\
&\qquad\qquad\qquad\qquad 
\text{(walks $\Zbar^u_k$ and $\Zbar^v_k$ are
independent under a common $\om$)}\\
&\leq C\Ef \biggl[ \bigl(E^\om g(T_{\Zbar_{k_1}}\om)  
\bigr)^2  \bigl(E^\om g(T_{\Zbar_{k_2}}\om)  \bigr)^2\\
&\qquad\qquad\qquad \times 
\sum_{u,v\in\Zb^2} P^\om\{\Zbar_{k_2+1}=u, \Zbar'_{k_2+1}=v\}
\biggr] P(Y^0_{k_3-k_2-1}=0)  \\
&=C\Ef \Bigl[ \bigl(E^\om g(T_{\Zbar_{k_1}}\om)  
\bigr)^2  \bigl(E^\om g(T_{\Zbar_{k_2}}\om)  \bigr)^2
 \Bigr] P(Y^0_{k_3-k_2-1}=0). 
\end{align*}
Now repeat the same step, and ultimately arrive at 
\begin{align*}
&\sum_{0\leq k_1<k_2<k_3<n} \Ef \Bigl[ \bigl(E^\om g(T_{\Zbar_{k_1}}\om)  
\bigr)^2  \bigl(E^\om g(T_{\Zbar_{k_2}}\om)  \bigr)^2
 \bigl(E^\om g(T_{\Zbar_{k_3}}\om)  \bigr)^2
\Bigr]\\
&\quad\leq C\sum_{0\leq k_1<k_2<k_3<n} P(Y^0_{k_1}=0) 
P(Y^0_{k_2-k_1-1}=0) 
P(Y^0_{k_3-k_2-1}=0)\\
&\quad\leq C\cdot G_{n-1}(0,0)^3 \leq Cn^{3/2}. \qedhere
\end{align*}
\end{proof}

By Theorem 8.8 in \cite[Chapter 3]{ethi-kurt},
\[
\Ef\bigl[\bigl( a_n(t+h,r)-a_n(t,r)\bigr)^3
\bigl( a_n(t,r)-a_n(t-h,r)\bigr)^3\,\bigr]\leq 
Ch^{3/2}
\]
is sufficient for  tightness of the 
processes $\{a_n(t,r):t\geq 0\}$.   The left-hand side
above is bounded by 
\[
\Ef\bigl[\bigl( a_n(t+h,r)-a_n(t,r)\bigr)^6\,\bigr]
+
\Ef\bigl[\bigl( a_n(t,r)-a_n(t-h,r)\bigr)^6\,\bigr].
\]
Note that if $h<1/(2n)$ then
 $\bigl( a_n(t+h,r)-a_n(t,r)\bigr)
\bigl( a_n(t,r)-a_n(t-h,r)\bigr)=0$ due to the discrete time
of the unscaled walks,  while if $h\geq 1/(2n)$ 
then $\tfl{n(t+h)}-\tfl{nt}\leq 3nh$. Putting these points
together shows that tightness will follow from the next moment 
bound.

\begin{lemma}
 There exists a constant $0<C<\infty$ such that,
for all $0\leq m <n\in\Nb$,  
\[
\Ef\bigl[ \bigl( \{E^\om(Z^{0,0}_n)-nV\}
- \{E^\om(Z^{0,0}_m)-mV\}  \bigr)^6\,\bigr]\leq C(n-m)^{3/2}. 
\]
\label{qm-mom-6-1}
\end{lemma} 

\begin{proof} The claim reduces  to Lemma \ref{qm-mom-6}
by restarting the walks at time $m$. 
\end{proof}

Convergence of finite-dimensional distributions 
in Corollary \ref{a-cor-1} follows from Theorem \ref{an-thm-1}.  
The limiting process $\abar(\cdot)=\lim a_n(\cdot,r)$ 
is identified by its covariance
$E\abar(s)\abar(t)=\Gamma_q\bigl((s\land t,r),(s\land t,r)\bigr)$. 
This completes the  proof of  Corollary \ref{a-cor-1}.  

\def\xn{X^{x(n,r)+\tfl{ntb},\,\tfl{nt}}_{\tfl{nt}}}
\def\xnsq{X^{x(n,q)+\tfl{nsb},\,\tfl{ns}}_{\tfl{ns}}}
\def\xns{X^{x(n,r)+\tfl{nsb},\,\tfl{ns}}_{\tfl{ns}}}
\def\xnk{X^{x(n,r_k)+\tfl{nt_kb},\,\tfl{nt_k}}_{\tfl{nt_k}}}
\def\xnl{X^{x(n,r_l)+\tfl{nt_lb},\,\tfl{nt_l}}_{\tfl{nt_l}}}

\section{Proofs for the random average process}
\label{pf-rap-sec} 

This section requires Theorem \ref{yn-thm-1}  from the space-time RWRE
section. 

\subsection{Separation of effects}
As the form of the limiting process in Theorem \ref{rap-thm-1}
 suggests, 
we can separate the fluctuations that come from the initial
configuration from those created by the dynamics.  
The quenched means of the RWRE represent the latter. We
start with the appropriate decomposition.  
Abbreviate 
\[
x_{n,r}=x(n,r)=\tfl{n\ybar}+\tfl{r\sqrt{n}\,}.
\]
Recall that we are considering $\ybar\in\Rb$ fixed, 
while $(t,r)\in\Rb_+\times\Rb$
is variable and serves as the index for the process.  
\begin{align*}
&\si_{\tfl{nt}}^n(x_{n,r}+\tfl{ntb})-\si_0^n(x_{n,r})=
E^\om\Bigl[\si_0^n(\xn)-\si_0^n(x_{n,r})\Bigr]\\
&=E^\om\biggl[{\bf1}_{\left\{\xn>x(n,r)\right\}}\sum_{i=x(n,r)+1}^{\xn}\eta_0^n(i)\\
&\qquad\qquad\qquad\qquad
 -{\bf1}_{\left\{\xn<x(n,r)\right\}}\sum_{i=\xn+1}^{x(n,r)}\eta_0^n(i)\biggr]\\
 &=\sum_{i>x(n,r)}P^\om\left\{i\leq\xn\right\}\cdot\eta_0^n(i)-\sum_{i\leq
x(n,r)}P^\om\left\{i>\xn\right\}\cdot\eta_0^n(i). 
 \end{align*}
 Recalling the means   $\vr(i/n)=\Ev\eta_0^n(i)$ we write this as
\begin{equation}
\si_{\tfl{nt}}^n(x_{n,r}+\tfl{ntb})-\si_0^n(x_{n,r})=
Y^n(t,r)+H^n(t,r)
 \label{eq:hy}
\end{equation}
 where 
\[
\begin{split}
Y^n(t,r)= \sum_{i\in\Zb}\bigl(\eta_0^n(i)-\vr(i/n)\bigr)
&\Bigl( \ind\{i>x_{n,r}\} P^\om\bigl\{ i\leq\xn\bigr\}  \\
&\qquad  -\; \ind\{i\leq x_{n,r}\} P^\om\bigl\{ i>\xn\bigr\}
\Bigr) 
\end{split}
\] 
and 
\[
\begin{split}
H^n(t,r)= \sum_{i\in\Zb} \vr(i/n) 
&\Bigl( \ind\{i>x_{n,r}\} P^\om\bigl\{ i\leq\xn\bigr\}  \\
&\qquad  -\; \ind\{i\leq x_{n,r}\} P^\om\bigl\{ i>\xn\bigr\}
\Bigr).  
\end{split}
\]

The plan of the proof of Theorem \ref{rap-thm-1}
is summarized in the next lemma. 
In the pages that follow we then show  the 
finite-dimensional weak convergence
$n^{-1/4}H^n \to H $, and   the finite-dimensional weak convergence
$n^{-1/4}Y^n \to Y $
 for a fixed $\om$.   This last statement is actually
not proved quite in the  strength just stated, but
the spirit is correct. 
The distributional limit $n^{-1/4}Y^n \to Y $ comes from 
the centered  initial increments  $\eta_0^n(i)-\vr(i/n)$, while 
a homogenization effect takes place for the coefficients
$P^\om\{i\leq\xn\}$ which converge to limiting deterministic
Gaussian probabilities. 
 Since the initial height functions 
$\si^n_0$ and the random environments $\om$ that drive the 
dynamics are independent, we also get  convergence 
$n^{-1/4}(Y^n+H^n)\to Y+H$ 
with  independent terms $Y$ and $H$. This is exactly
the statement of  Theorem \ref{rap-thm-1}.

\begin{lemma}
Let $(\Omega_0,\Fc_0,P_0)$ be a probability space
on which are defined independent random variables $\eta$ and
$\omega$ with values in some abstract measurable spaces.
The marginal laws are  $\Pf$  for $\omega$
and $\Pv$ for $\eta$, and   $\Pv^\om=\delta_\om\otimes\Pv$ 
is  the   conditional
probability distribution of $(\om,\eta)$, given $\om$. 
Let $\Hvec^n(\omega)$  and $\Yvec^n(\omega,\eta)$
be $\Rb^N$-valued measurable functions of $(\omega,\eta)$.
Make assumptions {\rm (i)--(ii)} below. 
\begin{itemize}
\item[(i)] There exists an $\Rb^N$-valued random vector
 $\Hvec$ such that
  $\Hvec^n(\omega)$  converges weakly  to 
 $\Hvec$.
\item[(ii)] There exists an $\Rb^N$-valued random vector
  $\Yvec$ such that, for all $\theta\in\Rb^N$,  
\[
\Ev^\om[e^{i\theta\cdot\Yvec^n}]\to 
E(e^{i\theta\cdot\Yvec}) \quad\text{ in $\Pf$-probability 
as $n\to\infty$.}
\]
 \end{itemize}
Then  
$\Hvec^n+\Yvec^n$ converges weakly to  
$\Hvec+\Yvec$, 
where $\Hvec$ and $\Yvec$ are independent.
\label{main-lm}
\end{lemma}

\begin{proof} 
 Let $\theta,\lambda$ be arbitrary vectors in $\bbR^N$. Then
\begin{align*}
&\left\vert \Ef\Ev^\om[e^{i\lambda\cdot\Hvec^n+i\theta\cdot\Yvec^n}]
- E[e^{i\lambda\cdot\Hvec}]\, E[e^{i\theta\cdot\Yvec}]\right\vert\\
&\qquad\leq
\left\vert \Ef\left[e^{i\lambda\cdot\Hvec^n}
\left(\Ev^\omega e^{i\theta\cdot\Yvec^n} - E e^{i\theta\cdot\Yvec}\right)\right]
\right\vert  +
\left\vert 
\left(\Ef e^{i\lambda\cdot\Hvec^n} - E e^{i\lambda\cdot\Hvec}\right)
E e^{i\theta\cdot\Yvec}
\right\vert\\
&\qquad\leq
\left\vert \Ef\left[e^{i\lambda\cdot\Hvec^n}
\left(\Ev^\omega e^{i\theta\cdot\Yvec^n} - E e^{i\theta\cdot\Yvec}\right)\right]
\right\vert  +
\left\vert
\Ef e^{i\lambda\cdot\Hvec^n} - E e^{i\lambda\cdot\Hvec}
\right\vert.
\end{align*}
By assumption (i), the second term above goes to 0. By assumption (ii), the 
integrand in the first term goes to 0 in $\Pf$-probability. Therefore by
bounded convergence the first term goes to 0 as $n\to\infty$.
\end{proof}

Turning to the work itself,  we check first that   $H^n(t,r)$ 
can be replaced 
with a quenched RWRE mean. 
Then the convergence  $H^n\to H$ follows from the RWRE results. 
 
\begin{lemma}
For any $S, T<\infty$ and for $\Pf$-almost every $\om$,
\[ 
\lim_{n\to\infty} \sup_{\substack{0\leq t\leq T\\-S\leq r\leq S}} 
n^{-1/4} \Bigl\lvert  
H^n(t,r)-\vr(\bar y)\cdot E^\om\bigl(\xn-x_{n,r}\bigr)
 \Bigr\rvert =0. 
\]
\label{H-rwre-lm}
\end{lemma}
 
\begin{proof} Decompose  $H^n(t,r)=H^n_1(t,r)-H^n_2(t,r)$ where 
\begin{align*}
 H^n_1(t,r)&=\sum_{i>x(n,r)}P^\om\bigl\{i\leq\xn\bigr\}\cdot
\vr(i/n),\\
H^n_2(t,r)&=\sum_{i\leq x(n,r)}P^\om\bigl\{i>\xn\bigr\}\cdot
\vr(i/n).
\end{align*}
  Working with $H^n_1(t,r)$, we separate out the negligible error. 
\begin{align*}
H^n_1(t,r)
&=\quad\vr(\bar y)\sum_{i>x(n,r)}P^\om\left\{i\leq\xn\right\}\\
&\qquad\qquad
+\sum_{i>x(n,r)}P^\om\left\{i\leq\xn\right\}\cdot\left[\vr(i/n)-\vr(\bar y)\right]\\
&=\quad\vr(\bar y)\cdot E^\om\left[\left(\xn-x_{n,r}\right)^+\right]
+R_1(t,r)
\end{align*}
with
\[
R_1(t,r)=\sum_{m=1}^\infty P^\om\left\{x_{n,r}+m\leq\xn\right\}\cdot\left[\vr\left(\frac{x_{n,r}}{n}+\frac{m}{n}\right)-\vr(\bar y)\right].
\]
Fix a small positive number $\de<\tfrac12$, and   use the  boundedness of
probabilities and the function $\vr$. 
\begin{multline}
|R_1(t,r)|\leq\sum_{m=1}^{\tfl{n^{1/2+\de}}}\left|\vr\left(\frac{x_{n,r}}{n}+\frac{m}{n}\right)-\vr(\bar y)\right|\\
+C\cdot\sum_{m=\tfl{n^{1/2+\de}}+1}^\infty P^\om\left\{x_{n,r}+m\leq\xn\right\}.\label{eq:r1}
\end{multline}
By the local H\"older-continuity  of $\vr$ with exponent $\gamma>\frac12$, the first sum is $\mathfrak{o}(n^{1/4})$ 
 if $\de>0$ is small enough. 
Since $X^{x(n,r)+\tfl{ntb}, \tfl{nt}}_0=x_{n,r}+\tfl{ntb}$ 
and by time $\tfl{nt}$ the walk has displaced by at most $M\tfl{nt}$,
there are at most $\Oc(n)$ nonzero terms in 
the second sum in \eqref{eq:r1}. Consequently this sum is at most  
\[   Cn\cdot P^\om\left\{\xn-x_{n,r}\geq \tfl{n^{1/2+\de}}\right\}.  
\]
  By Lemma \ref{lm:moddev}  the last line vanishes  
uniformly over  $t\in[0,T]$ and $r\in[-S,S]$ 
 as $n\to\infty$,
 for $\Pf$-almost every $\om$.  
We have shown 
\[ 
\lim_{n\to\infty} \sup_{\substack{0\leq t\leq T\\-S\leq r\leq S}} n^{-1/4} \Bigl\lvert  
H^n_1(t,r) -\vr(\bar y)\cdot E^\om\bigl[\bigl(\xn-x_{n,r}\bigr)^+\,\bigr]
 \Bigr\rvert =0 \quad\text{$\Pf$-a.s.} 
\]
Similarly one shows 
\[
\lim_{n\to\infty} \sup_{\substack{0\leq t\leq T\\-S\leq r\leq S}} n^{-1/4} \Bigl\lvert  
H^n_2(t,r) -\vr(\bar y)\cdot E^\om\bigl[\bigl(\xn-x_{n,r}\bigr)^-\,\bigr]
 \Bigr\rvert =0   \quad\text{$\Pf$-a.s.}
 \]
The conclusion follows from the combination of these two. 
 \end{proof}

For a fixed $n$ and $\ybar$, the process 
$E^\om\bigl(\xn-x_{n,r}\bigr)$ has the 
same distribution  as the process $y_n(t,r)$ defined in 
\eqref{def-yn}.  A combination of Lemma \ref{H-rwre-lm}
and Theorem \ref{yn-thm-1} imply that 
the finite-dimensional distributions of the processes
$n^{-1/4}H_n$ converge weakly, as $n\to\infty$,
to the  finite-dimensional distributions of the mean-zero
Gaussian  process $H$ with covariance 
\be
EH(s,q)H(t,r)
=\vr(\ybar)^2 \Gamma_q((s,q),(t,r)). 
\label{H-covar}
\ee

\subsection{Finite-dimensional convergence of $Y^n$}

Next we turn to  convergence of the 
finite-dimensional distributions of process $Y^n$  in \eqref{eq:hy}.
 Recall 
that   $B(t)$ is  standard Brownian motion, and 
$\sigma_a^2=E[(X^{0,0}_1-V)^2]$ is
the variance
of the annealed walk.  Recall the definition 
\[
\begin{split}
&\Gamma_0((s,q),(t,r))
=\int_{q\lor r}^\infty 
P[ \sigma_a  B(s)> x-q]P[  \sigma_a B(t)> x-r]\,dx \\
&\qquad\qquad 
-  \Bigl\{ \ind_{\{r>q\}} \int_q^r P[   \sigma_a B(s)> x-q]
P[  \sigma_a B(t)\leq x-r]\,dx \\
&\qquad\qquad\qquad\qquad + 
\ind_{\{q>r\}} \int_r^q P[   \sigma_a B(s)\leq x-q]
P[  \sigma_a B(t)> x-r]\,dx \Bigr\} \\
&\qquad\qquad + 
 \int_{-\infty}^{q\land r} P[   \sigma_a B(s)\leq x-q]
P[  \sigma_a B(t)\leq x-r]\,dx. 
\end{split}
\]
Recall from \eqref{eq:meta}
 that $v(\ybar)$ is the variance of the increments
around $\tfl{n\ybar}$.  
Let $\{Y(t,r):t\geq 0,\, r\in\Rb\}$ be a real-valued  mean-zero Gaussian 
process with covariance
\be
EY(s,q)Y(t,r)=v(\ybar)\Gamma_0((s,q),(t,r)).
\label{Y-covar}
\ee
  Fix $N$ 
and  space-time points $(t_1,r_1),\dotsc,(t_N,r_N)\in\Rb_+\times\Rb$.
Define vectors 
\[
\Yvec^n=n^{-1/4}\bigl(Y^n(t_1,r_1),\dotsc, Y^n(t_N,r_N)\bigr) 
\quad\text{and}\quad
\Yvec= \bigl(Y(t_1,r_1),\dotsc, Y(t_N,r_N)\bigr).  
\]
This section is devoted to the proof of the next proposition,
after which we finish the proof of Theorem \ref{rap-thm-1}. 

\begin{proposition}
For any vector $\theta=(\theta_1,\dotsc,\theta_N)\in\Rb^N$,
$\Ev^\w(e^{i\theta\cdot\Yvec^n})\to 
 E(e^{i \theta\cdot\Yvec })$  in $\Pf$-probability
as $n\to\infty$.
\label{Y-findim-prop}
\end{proposition}

\begin{proof} 
Let $G$ be a centered 
Gaussian variable with 
variance
\[
S=v(\ybar)
 \sum_{k,\,l=1}^N\te_k\te_l \Gamma_0((t_k,r_k),(t_l,r_l))
\]
 and so $\theta\cdot\Yvec$ is distributed like 
$G$. 
We will show that 
\be
\Ev^\w(e^{i\theta\cdot\Yvec^n  })\to 
 E(e^{i G}) \quad\text{ in $\Pf$-probability.}  
\label{goal-Yfd-1}
\ee

Recalling the definition of $Y^n(t,r)$,
introduce some notation: 
\[
\begin{split}
\zeta^\om_n(i,t,r)&=\ind\{i>x_{n,r}\} P^\om\bigl\{ i\leq\xn\bigr\}\\  
&\qquad\qquad 
\;-\; \ind\{i\leq x_{n,r}\} P^\om\bigl\{ i>\xn\bigr\} 
\end{split}
\]
so that 
\[
Y^n(t,r)= \sum_{i\in\Zb}\bigl(\eta_0^n(i)-\vr(i/n)\bigr)
\zeta^\om_n(i,t,r). 
\]
Then put 
\[
\nu^\om_n(i)=\sum_{k=1}^N\te_k \,\zeta^\om_n(i,t_k,r_k) 
\]
and
\[
U_n(i)=n^{-1/4}\left(\eta_0^n(i)-\vr(i/n)\right)\nu^\om_n(i).
\]
Consequently 
\[\theta\cdot\Yvec^n=\sum_{i\in\Zb } U_n(i).\]
To separate out the relevant terms let $\delta>0$ be small 
and  define 
\[
W_n=\sum_{i=\tfl{n\ybar}-\tfl{n^{1/2+\de}}}^{\tfl{n\ybar}+\tfl{n^{1/2+\de}}}
 U_n(i).\]
For fixed $\om$ and $n$,  under the measure $\Pv^\om$ 
the variables $\{U_n(i)\}$ are constant multiples of
centered increments $\eta_0^n(i)-\vr(i/n)$ and hence 
independent and mean zero. Recall also that second moments of
centered increments $\eta_0^n(i)-\vr(i/n)$ are uniformly bounded. 
Thus the   terms left out of $W_n$ satisfy 
\begin{align*}
\Ev^\om\bigl[(W_n- \theta\cdot\Yvec^n)^2\,\bigr]
\leq  Cn^{-1/2}\sum_{i: \lvert i- \tfl{n\ybar}\,\rvert
\,>\,  {n^{1/2+\de}}} 
\nu^\om_n(i)^2,
\end{align*}
and we wish to show that this upper bound vanishes 
 for $\Pf$-almost every $\om$ as $n\to\infty$. Using the definition
of $\nu^\om_n(i)$, bounding the sum on the right
 reduces to bounding sums of the two types
\[
n^{-1/2}\sum_{i: \lvert i- \tfl{n\ybar}\,\rvert
\,>\,  {n^{1/2+\de}}} \ind\{i>x(n,r_k)\} 
\Bigl( P^\om\bigl\{ i\leq\xnk\bigr\}\Bigr)^2
\]
and 
\[
n^{-1/2}\sum_{i: \lvert i- \tfl{n\ybar}\,\rvert
\,>\,  {n^{1/2+\de}}}\ind\{i\leq x(n,r_k)\} 
\Bigl(P^\om\bigl\{ i>\xnk\bigr\}\Bigr)^2. 
\]
For large enough $n$ the points $x(n,r_k)$ lie within $\tfrac12 n^{1/2+\de}$
of $\tfl{n\ybar}$, and then the previous sums are bounded
by the sums 
\[
n^{-1/2}\sum_{i\,\geq\,  x(n,r_k) + (1/2)  {n^{1/2+\de}}}  
\Bigl( P^\om\bigl\{ i\leq\xnk\bigr\}\Bigr)^2
\]
and 
\[
n^{-1/2}\sum_{i\,\leq\,  x(n,r_k)  -  (1/2) {n^{1/2+\de}}} 
\Bigl(P^\om\bigl\{ i>\xnk\bigr\}\Bigr)^2. 
\]
These vanish  for $\Pf$-almost every $\om$ as $n\to\infty$ by
Lemma \ref{lm:moddev}, 
in a manner similar to  the second sum in \eqref{eq:r1}.  
Thus $\Ev^\om\bigl[(W_n- \theta\cdot\Yvec^n)^2\,\bigr]\to 0$
and our goal \eqref{goal-Yfd-1} has simplified to 
\be
\Ev^\w(e^{i W_{n}})\to 
 E(e^{i G}) \quad\text{ in $\Pf$-probability.}  
\label{goal-Yfd-2}
\ee

We use  the Lindeberg-Feller theorem to formulate 
conditions for  a central 
limit theorem for $W_n$ under a fixed $\om$. 
For Lindeberg-Feller we need to check two conditions:
\begin{itemize}
\item[(LF-i)]
$\ds{\; S^n(\om)\;\equiv 
\sum_{i=\tfl{n\ybar}-\tfl{n^{1/2+\de}}}^{\tfl{n\ybar}+\tfl{n^{1/2+\de}}}
\Ev^\om\left[U_n(i)^2\right]
\ds\mathop{\longrightarrow}_{n\to\infty}S} $\\
\item[(LF-ii)]
$\ds{\sum_{i=\tfl{n\ybar}-\tfl{n^{1/2+\de}}}^{\tfl{n\ybar}+\tfl{n^{1/2+\de}}}
\Ev^\om\left[U_n(i)^2\cdot{\bf 1}_{\{|U_n(i)|>\ve\}}\right]
\ds\mathop{\longrightarrow}_{n\to\infty}0}\;$
for all $\ve>0$.
\end{itemize}
To  see that (LF-ii) holds, 
pick conjugate exponents $p, q>1$  ($1/p+1/q=1$):
\begin{align*}
&\Ev^\om\left[U_n(i)^2\cdot{\bf 1}_{\{U_n(i)^2>\ve^2\}}\right]
\leq 
\left(\Ev^\om\left[|U_n(i)|^{2p}\right]\right)^\frac{1}{p}
\left(\Pv^\om\left[U_n(i)^2>\ve^2\right]\right)^\frac{1}{q}\\
&\leq \ve^{-\frac2q} \left(\Ev^\om\left[|U_n(i)|^{2p}\right]
\right)^\frac{1}{p}
\left(\Ev^\om\left[U_n(i)^2\right]\right)^\frac{1}{q}
\leq Cn^{-1/2-1/(2q)}.
\end{align*}
In the last step we used the bound $|U_n(i)|\leq 
Cn^{-1/4}\left|\eta_0^n(i)-\vr(i/n)\right|$, boundedness
of $\vr$, 
and we took  $p$ close enough to $1$ to apply assumption
\eqref{eta-ass-p}.  
Condition (LF-ii)  follows if
$\de<1/(2q)$.

We turn to condition (LF-i).  
\begin{align*}
&S^n(\om)=
\sum_{i=\tfl{n\ybar}-\tfl{n^{1/2+\de}}}^{\tfl{n\ybar}+\tfl{n^{1/2+\de}}}
\Ev^\om\left[U_n(i)^2\right]
=\sum_{i=\tfl{n\ybar}-\tfl{n^{1/2+\de}}}^{\tfl{n\ybar}+\tfl{n^{1/2+\de}}}n^{-1/2}v
(i/n)[\nu^\om_n(i)]^2\\
&=
\sum_{i=\tfl{n\ybar}-\tfl{n^{1/2+\de}}}^{\tfl{n\ybar}+\tfl{n^{1/2+\de}}}n^{-1/2}
\left[v(i/n)-v(\bar y)\right][\nu^\om_n(i)]^2
+\sum_{i=\tfl{n\ybar}-\tfl{n^{1/2+\de}}}^{\tfl{n\ybar}+\tfl{n^{1/2+\de}}}n^{-1/2}v(\bar y)[\nu^\om_n(i)]^2
\end{align*}
Due to the local H\"older-property \eqref{holder} of $v$, the first sum
on the last line is bounded above by
\[
C(\bar y)n^{1/2+\de}n^{-1/2}\left[n^{-1/2+\de}\right]^\gamma
=C(\bar y)n^{\de(1+\gamma)-\gamma/2}\to0
\]
for sufficiently small $\de$. Denote the   remaining  relevant
part by $\tilde{S}^n(\om)$, given by 
\begin{align}
&\tilde S^n(\w)=
\sum_{i=\tfl{n\ybar}-\tfl{n^{1/2+\de}}}^{\tfl{n\ybar}
+\tfl{n^{1/2+\de}}} n^{-1/2}v(\bar y)[\nu^\om_n(i)]^2
=v(\bar y)n^{-1/2}\sum_{m=-\tfl{n^{1/2+\de}}}^{\tfl{n^{1/2+\de}}}
\left(\nu^\om_n\left(m+\tfl{n\ybar}\right)\right)^2\nn\\
=&\,v(\bar y)\sum_{k,\,l=1}^N\te_k\te_l\; n^{-1/2}
\sum_{m=-\tfl{n^{1/2+\de}}}^{\tfl{n^{1/2+\de}}}
\zeta^\om_n(\tfl{n\ybar}+m,t_k,r_k)\zeta^\om_n(\tfl{n\ybar}+m,t_l,r_l).  
\label{def-Sn-Y}
\end{align}
Consider for the moment a particular $(k,l)$ term 
in the  first sum on line  \eqref{def-Sn-Y}. Rename
$(s,q)=(t_k,r_k)$ and $(t,r)=(t_l,r_l)$. 
Expanding the product of the $\zeta^\om_n$-factors gives 
three sums: 
\begin{align}
&n^{-1/2}
\sum_{m=-\tfl{n^{1/2+\de}}}^{\tfl{n^{1/2+\de}}}
\zeta^\om_n(\tfl{n\ybar}+m,s,q)\zeta^\om_n(\tfl{n\ybar}+m,t ,r )\nn\\
&=n^{-1/2}
\sum_{m=-\tfl{n^{1/2+\de}}}^{\tfl{n^{1/2+\de}}}
\ind_{\{m>\tfl{q\sqrt{n}\,}\}}\ind_{\{m>\tfl{r\sqrt{n}\,}\}}
P^\om\bigl( \xnsq \geq \tfl{n\ybar}+m\bigr) \nn\\
&\qquad \qquad\qquad \qquad\qquad \qquad
 \times P^\om\bigl( \xn \geq \tfl{n\ybar}+m\bigr)
\label{sum-S-1}\\
&-\; n^{-1/2} \sum_{m=-\tfl{n^{1/2+\de}}}^{\tfl{n^{1/2+\de}}} \biggl\{
\ind_{\{m>\tfl{q\sqrt{n}\,}\}}\ind_{\{m\leq \tfl{r\sqrt{n}\,}\}}
P^\om\bigl( \xnsq \geq \tfl{n\ybar}+m\bigr) \nn\\
&\qquad \qquad\qquad \qquad\qquad \qquad
 \times P^\om\bigl( \xn < \tfl{n\ybar}+m\bigr)\nn\\
&\qquad\qquad\qquad
 + \;\ind_{\{m\leq \tfl{q\sqrt{n}\,}\}}\ind_{\{m> \tfl{r\sqrt{n}\,}\}}
P^\om\bigl( \xnsq < \tfl{n\ybar}+m\bigr) \nn\\
&\qquad \qquad\qquad \qquad\qquad \qquad
 \times P^\om\bigl( \xn \geq \tfl{n\ybar}+m\bigr)
\biggr\}  \label{sum-S-2}\\
&+\;n^{-1/2}
\sum_{m=-\tfl{n^{1/2+\de}}}^{\tfl{n^{1/2+\de}}}
\ind_{\{m\leq \tfl{q\sqrt{n}\,}\}}\ind_{\{m\leq\tfl{r\sqrt{n}\,}\}}
P^\om\bigl( \xnsq < \tfl{n\ybar}+m\bigr) \nn\\
&\qquad \qquad\qquad \qquad\qquad \qquad
 \times P^\om\bigl( \xn < \tfl{n\ybar}+m\bigr)
\label{sum-S-3}
\end{align}
Each of these three sums
\eqref{sum-S-1}--\eqref{sum-S-3} converges to a corresponding integral
in $\Pf$-probability, due to the quenched CLT Theorem \ref{q-inv-pr}.
To see the correct limit, just note that 
\[
\begin{split}
&P^\om\bigl( \xn < \tfl{n\ybar}+m\bigr)\\
&\qquad\qquad =P^\om\bigl( \xn -X^{x(n,r)+\tfl{ntb},\tfl{nt}}_0
 < -\tfl{ntb}+m-\tfl{r\sqrt{n}\,}\;\bigr)\end{split}
\]
and recall that $-b=V$ is the average speed of the walks.
 We give technical details of  the argument for the first sum
in the next lemma.

\begin{lemma} As $n\to\infty$, 
the sum in {\rm\eqref{sum-S-1}} converges
in $\Pf$-probability to
\[ 
\int_{q\lor r}^\infty
P[ \sigma_a B(s)> x-q]P[ \sigma_a B(t)> x-r]\,dx.
\]
\label{S-conv-lm}
\end{lemma}

\begin{proof}[Proof of Lemma \ref{S-conv-lm}]
With 
\begin{align*}
f_n^\w(x)=
P^\w\bigl(  \xnsq\geq \tfl{n\ybar}+\tfl{x\sqrt n}\bigr)
P^\w\bigl(  \xn\geq \tfl{n\ybar}+\tfl{x\sqrt n}\bigr) 
\end{align*}
and 
 \[I^\om_n=\int_{q\lor r}^{n^\delta} f_n^\w(x)dx,\]
the sum in \eqref{sum-S-1} equals $I^\om_n+{\mathcal O}(n^{-1/2})$.
 By the quenched invariance principle 
 Theorem \ref{q-inv-pr}, 
for any   fixed $x$, $f_n^\w(x)$
converges in $\Pf$-probability to
\[
 f(x)=P[ \sigma_a B(s)\geq x-q]P[ \sigma_a B(t)\geq x-r].\]
We cannot claim this  convergence $\Pf$-almost surely because
the walks $\xn$ change as $n$ changes. 
But by  a textbook characterization of convergence in probability,
 for  a fixed $x$  each subsequence $n(j)$ 
has a further subsequence $n(j_\ell)$ such that
\[
\Pf\Bigl[\om\,:\;
f_{n(j_\ell)}^\w(x)\underset{\ell\to\infty}{\longrightarrow}
f(x)\Bigr]=1.
\]
By the diagonal trick, one can find one subsequence for all $x\in\bbQ$ 
and thus 
\[\forall\{n(j)\}, \exists\{j_\ell\}:
\Pf\Bigl[\om\,:\;\forall x\in\bbQ:f_{n(j_\ell)}^\w(x)\to f(x)\Bigr]=1.\]
Since $f_n^\w$ and $f$ are nonnegative and nonincreasing, and
 $f$ is 
continuous and decreases to 0, 
 the convergence works for all $x$ and is uniform on $[q\lor r,\infty)$. 
That is,
\begin{align*}
\forall\{n(j)\}, \exists\{j_\ell\}:
\Pf\Bigl[\om\,:\;\norm{f_{n(j_\ell)}^\w-f}_{L^\infty[q\lor r,\,\infty)}
\to0\Bigr]=1.
\end{align*}
It remains to make the step to the convergence of the integral 
$I^\om_n$ to $\int_{q\lor r}^\infty f(x)\,dx$. 

Define now \[J^\om_n(A)=\int_{q\lor r}^A f^\w_n(x)dx.\] Then, for any 
 $A<\infty$
\begin{align*}
\forall\{n(j)\}, \exists\{j_\ell\}:
\Pf\Bigl[\om\,:\; J^\om_{n(j_\ell)}(A)\to \int_{q\lor r}^A f(x)dx\Bigr]=1.
\end{align*}
In other words, $J^\om_n(A)$ converges to $\int_{q\lor r}^A f(x)dx$ in $\Pf$-probability.
Thus, for each $0<A<\infty$, there is an 
integer $m(A)$ such that for all $n\geq m(A)$
\[\Pf\bigg[\om\,:\;
\Big|J^\om_n(A)-\int_{q\lor r}^A f(x)dx\Big|>A^{-1}\bigg]<A^{-1}.\]
Pick $A_n\nearrow\infty$ such that 
 $m(A_n)\leq n$.
Under the annealed measure
 $P$, $X_n^{0,0}$ is a homogeneous mean zero random walk
with variance $\Oc(n)$.  Consequently
\begin{align*}
&\Ef[\,\lvert I^\om_n-J^\om_n(A_n)\rvert\,]\leq
 \int_{A_n\land n^\delta}^\infty\Ef[f^\w_n(x)]dx
\nn\\
&\qquad \leq \int_{A_n\land n^\delta}^\infty
P\Bigl( \xn \geq x(n,r)-\tfl{r\sqrt n}+\tfl{x\sqrt n} \Bigr)\,dx
\mathop{\longrightarrow}_{n\to\infty}0.
 \end{align*}
Combine this with 
\[\Pf\Bigl[\om\,:\;
\Big|J^\om_n(A_n)-\int_{q\lor r}^{A_n} f(x)dx\Big|>A_n^{-1}\Bigr]<A_n^{-1}.\]
Since $\int_{q\lor r}^{A_n} f(x)dx$ converges to $\int_{q\lor r}^\infty f(x)dx$, we have shown that
$I^\om_n$ converges to this same integral in $\Pf$-probability. 
This completes the proof of Lemma \ref{S-conv-lm}. 
\end{proof} 

We return to the main development, the proof of Proposition  
\ref{Y-findim-prop}.  Apply the argument of the lemma to the 
three sums \eqref{sum-S-1}--\eqref{sum-S-3} to conclude the
following  
limit in $\Pf$-probability. 
\begin{align*}
& \lim_{n\to\infty} \;  n^{-1/2}
\sum_{m=-\tfl{n^{1/2+\de}}}^{\tfl{n^{1/2+\de}}}
\zeta^\om_n(\tfl{n\ybar}+m,s,q)\zeta^\om_n(\tfl{n\ybar}+m,t ,r )\\
&= \int_{q\lor r}^\infty 
P[ \sigma_a  B(s)> x-q]P[  \sigma_a B(t)> x-r]\,dx \\
&\qquad\qquad 
-  \Bigl\{ \ind_{\{r>q\}} \int_q^r P[   \sigma_a B(s)> x-q]
P[  \sigma_a B(t)\leq x-r]\,dx \\
&\qquad\qquad\qquad\qquad + 
\ind_{\{q>r\}} \int_r^q P[   \sigma_a B(s)\leq x-q]
P[  \sigma_a B(t)> x-r]\,dx \Bigr\} \\
&\qquad\qquad + 
 \int_{-\infty}^{q\land r} P[   \sigma_a B(s)\leq x-q]
P[  \sigma_a B(t)\leq x-r]\,dx\\
&=\Gamma_0((s,q),(t,r)). 
\end{align*}

Return to condition (LF-i) of the Lindeberg-Feller theorem
and the definition \eqref{def-Sn-Y} of $\tilde S^n(\om)$.
Since $S^n(\om)-\tilde S^n(\om)\to 0$ as pointed out above
\eqref{def-Sn-Y}, we have shown that $S^n\to S$ in
$\Pf$-probability. Consequently  
\begin{align*}
\forall\{n(j)\},\exists\{j_\ell\}:
\Pf\Bigl[\om\,:\;S^{n(j_\ell)}(\w)\to S\Bigr]=1.
\end{align*}
This can be rephrased as: given any subsequence
$\{n(j)\}$, there exists a further subsequence
 $\{n(j_\ell)\}$ 
along which conditions (LF-i) and (LF-ii) of the Lindeberg-Feller
theorem are satisfied for the array 
\[\{U_{n(j_\ell)}(i):\tfl{n(j_\ell)\ybar}-\tfl{n(j_\ell)^{1/2+\de}} \leq i
\leq \tfl{n(j_\ell)\ybar}+\tfl{n(j_\ell)^{1/2+\de}}\,,\; \ell\geq 1 \}\]
 under the measure $\Pv^\om$
 for $\Pf$-a.e.\ $\w$.
This implies that 
\[\forall\{n(j)\},\exists\{j_\ell\}:
\Pf\Bigl[\om\,:\;\Ev^\w(e^{i W_{n(j_\ell)}})\to 
 E(e^{i G})\Bigr]=1.\]
But the last statement
 characterizes convergence $\Ev^\w(e^{i W_{n}})\to 
 E(e^{i G})$ in $\Pf$-probability.  
As we already showed above that  $W_n- \theta\cdot\Yvec^n\to 0$
 in $\Pv^\om$-probability  $\Pf$-almost surely,   
this  completes the proof of Proposition \ref{Y-findim-prop}.
\end{proof}

\subsection{Proofs of  Theorem \ref{rap-thm-1} and Proposition
\ref{prop-2pt}}

\begin{proof}[Proof of Theorem \ref{rap-thm-1}]
The decomposition \eqref{eq:hy} gives 
$z_n=n^{-1/4}(Y^n+H^n)$. The paragraph that follows
Lemma \ref{H-rwre-lm} and  Proposition \ref{Y-findim-prop}
  verify the hypotheses
of Lemma \ref{main-lm} for $H^n$ and $Y^n$. 
Thus we have the limit $z_n\to z\equiv Y+H$ in the sense of 
convergence of finite-dimensional distributions.
Since 
$Y$ and $H$ are mutually  independent mean-zero Gaussian 
processes,  their covariances in 
\eqref{H-covar} and \eqref{Y-covar} can be added to give
\eqref{rap-cov}. 
\end{proof} 

\begin{proof}[Proof of Proposition \ref{prop-2pt}]
  The value  \eqref{2pt-beta} for $\beta$ 
 can be computed from \eqref{def-beta-1}, or from the 
probabilistic characterization \eqref{def-beta-2} of $\beta$
via  Example \ref{ex-abar-1}.  If we let $u$ denote a random
variable distributed like $u_0(0,-1)$, then we get 
\[
\beta= 
\frac{\Ef u-\Ef (u^2)}{\Ef u-(\Ef u)^2}
\quad\text{and}\quad
\kappa=\frac{\Ef (u^2)-(\Ef u)^2}{\Ef u-\Ef( u^2)}\,.
\]

With obvious notational simplifications, the evolution 
step \eqref{2pt-eta-evol} rewrites as 
\[
\eta'(k)-\rho= (1-u_{k}) (\eta(k)-\rho)+ u_{k-1} (\eta(k-1)-\rho)
+(u_{k-1}-u_{k})\rho. 
\]
Square both sides, take expectations, use the independence
of all variables $\{\eta(k-1),\eta(k), u_k, u_{k-1}\}$ on the right, 
and use the
requirement that $\eta'(k)$ have the same variance $v$ as 
$\eta(k)$ and $\eta(k-1)$. The result is the identity  
\[
v 
= v( 1-2\Ef u+2\Ef (u^2)) +2\rho^2 (\Ef (u^2)-(\Ef u)^2)
\] 
from which follows $v=\kappa \rho^2$.
The rest of part (b) is a straightforward specialization of 
\eqref{rap-cov}. 
\end{proof}










\bibliographystyle{abbrv} 
\bibliography{refs}
\end{document}